\newtheorem{lemma}{Lemma}[section]
\newtheorem{prop}[lemma]{Proposition}
\newtheorem{theorem}[lemma]{Theorem}
\newtheorem{cor}[lemma]{Corollary}
\newtheorem{rem}[lemma]{Remark}

\newcommand{\kla}{\left ( }
\newcommand{\mer}{\right ) }
\newcommand{\nach}{\rightarrow}
\newcommand{\nachop}[1]{\stackrel{#1}{\longrightarrow}}
\newcommand{\for}{\begin{eqnarray*}}
\newcommand{\mel}{\end{eqnarray*}}
\newcommand{\mitt}{\left | { \atop } \right.}
\newcommand{\kl}{\pl \le \pl}
\newcommand{\gl}{\pl \ge \pl}
\newcommand{\kll}{\p \le \p}

\newcommand{\summ}{\sum\limits}

\newcommand{\ew}{{\rm  I\! E}}
\newcommand{\nz}{{\rm  I\! N}}
\newcommand{\rz}{{\rm  I\! R}}
\newcommand{\kz}{{\rm  I\! K}}

\newcommand{\p}{\hspace{.05cm}}
\newcommand{\pl}{\hspace{.1cm}}
\newcommand{\pll}{\hspace{.3cm}}
\newcommand{\hz}{\vspace{0.5cm}}
\newcommand{\plll}{\hspace{1.1cm}}
\newcommand{\pla}{\hspace{1.5cm}}

\newcommand{\Om}{\Omega}
\newcommand{\om}{\omega}
\newcommand{\al}{\alpha}
\newcommand{\be}{\beta}
\newcommand{\La}{\Lambda}
\newcommand{\vare}{\varepsilon}

\newcommand{\lzn}{\ell_2^n}
\newcommand{\bzn}{B_{\lzn}}
\newcommand{\bfi}{B_{\Phi}}
\newcommand{\len}{\ell_1^n}
\newcommand{\lin}{\ell_{\infty}^n}
\newcommand{\efi}{E_{\Phi}}

\newcommand{\com}{{\cal L}(\lzn ,}
\newcommand{\comxy}{{\cal L}(X,Y)}

\newcommand{\izin}{\iota_{2,\infty}^n}
\newcommand{\iizn}{\iota_{\infty ,2}^n}
\newcommand{\iezn}{\iota_{1,2}^n}
\newcommand{\izen}{\iota_{2,1}^n}

\newcommand{\noo}{\left \|}
\newcommand{\rrm}{\right \|}
\newcommand{\bet}{\left |}
\newcommand{\rag}{\right |}
\newcommand{\ssum}{\sum\limits}

\newcommand{\pz}{\pi_2(}
\newcommand{\pzn}{\pi_2^n(}
\newcommand{\epz}{1\!\!1\!\!{\rm -}\! \pzn}
\newcommand{\fipsih}[1]{\delta (#1 |\Phi ,\Psi )}
\newcommand{\fty}[1]{t_{\Phi} (#1)}
\newcommand{\fco}[1]{c_{\Phi} (#1)}
\newcommand{\ftyh}[1]{\hat{t}_{\Phi} (#1)}

\newcommand{\pif}{\pi_{\Phi}}
\newcommand{\vf}{\nu_{\Phi}}

\newcommand{\on}{\Phi=(\phi_k)}
\newcommand{\fis}{\Phi = (\phi_k)_1^n}
\newcommand{\fiss}{\Phi = (\phi_k)_{k \in I}}

\newcommand{\psis}{\Psi = (\psi_k)_1^n}

\newcommand{\gmm}{\left\| \sum\limits_k g_k x_k \right\|_{L_2(X)}}

\newcommand{\Tgmm}{\left\| \sum\limits_k g_k T x_k \right\|_{L_2(Y)}}

\newcommand{\fmm}{\left\| \sum\limits_k \phi_k x_k \right\|_{L_2(X)}}

\newcommand{\Tfmm}{\left\| \sum\limits_k \phi_k T x_k \right\|_{L_2(Y)}}
\newcommand{\fme}{\left\| \sum\limits_{k=1}^n \phi_k x_k
\right\|_{L_1(\ell_1)}}

\documentstyle[12pt]{article}
\begin{document}


\title{Type and cotype with respect to arbitrary orthonormal systems}
\author{Stefan Geiss and Marius Junge
\thanks{The authors are supported by the DFG (Ko 962/3-1).}}
\date{}
\maketitle

{\center Mathematisches Institut der Friedrich Schiller Universit\H{a}t Jena\\
UHH 17. OG \\
D O7743 Jena \\
Germany \\
E-mail: geiss@mathematik.uni-jena.d400.de
\vspace{.7cm}

Mathematisches Seminar der Universit\H{a}t Kiel\\
Ludewig Meyn Strasse 4\\
D 24098 KIEL\\
Germany\\
E-mail: nms06@rz.uni-kiel.d400.de\\}



\begin{abstract}
Let $\on_{k \in \nz}$ be an orthonormal system on some $\sigma$-finite
measure space $(\Om,p)$. We study the notion of cotype with respect to $\Phi$
for an operator $T$ between two Banach spaces $X$ and $Y$, defined by
$\fco T := \inf$ $c$ such that
\[ \Tfmm \pl \le \pl c \pll \gmm \hspace{.7cm}\mbox{for all}\hspace{.7cm}
   (x_k)\subset X \pl,\]
where $(g_k)_{k\in \nz}$ is a sequence of independent and normalized
gaussian variables. It is shown that this $\Phi$-cotype coincides with
the usual notion of cotype $2$ iff \linebreak
$\fco {I_{\lin}} \sim \sqrt{\frac{n}{\log (n+1)}}$ uniformly in $n$
iff there is a positive $\eta>0$ such that for all
$n \in \nz$ one can find an orthonormal $\Psi = (\psi_l)_1^n \subset
{\rm span}\{ \phi_k \p|\p k \in \nz\}$ and a sequence of disjoint measurable
sets $(A_l)_1^n \subset \Om$ with
\[ \int\limits_{A_l} \bet \psi_l\rag^2 d p \pl \ge \pl \eta \quad
\mbox{for all}\quad l=1,...,n \pl. \]
A similar result holds for the type situation.
The study of type and cotype with respect to orthonormal systems
of a given length provides the appropriate approach to this result.
We intend to give a quite complete picture for orthonormal systems
in measure space with few atoms.
\end{abstract}

\newpage


\setcounter{section}{0}

\section*{Introduction and notation}
\newtheorem{ttheorem}{Theorem}
\newtheorem{example}[ttheorem]{Example}

The theory of type and cotype in Banach spaces is closely connected
to the probability theory and provides a good frame work to distinguish
relevant
local properties of Banach spaces. In this paper we develop a connection
between this theory and geometric properties of orthonormal systems. We do
this by means of a certain approximation of orthonormal systems by
systems of functions having disjoint supports. The main technical tool is the
use of operator ideal techniques.\hz

Throughout this paper the standard notation of the Banach space theory is
used. All Banach spaces are real or complex, in particular $\kz$ stands for
the real or complex scalars. Given a Banach space $X$ then $B_X$ denotes the
closed unit ball, $X^*$ the dual and $I_X$ the identity.
If $x\in X$ and $a\in X^*$ then $\langle x,a \rangle := a(x)$,
whereas for the scalar product in a Hilbert space we use $(\cdot ,\cdot )$.
If $(\Om ,p)$ is a $\sigma$-finite measure space (we will always
assume $\sigma$-finiteness) and if $1\le r<\infty$ then $L_r(X)$ is the
Banach space of all $X$-valued strongly measurable functions $f:\Om\nach X$
such that $\| f\|_{L_r(X)} := \kla \int_{\Om} \|f\|^r d p \mer^{1/r} <
\infty$. For two Banach spaces $X$ and $Y$ as usual $\comxy$ is
the Banach space of all linear and continuous operators from $X$ into $Y$
equipped with the operator norm
$\| Tx\| = \sup \left\{ \| Tx \| \mitt x\in B_X \right\}$.
To shorten some statements let us denote the formal identity
$\ell_p^n \nach \ell_q^n$ by $\iota_{pq}^n$.
In the whole paper $(g_k)_{k\in \nz}$
is a sequence of independent, normalized gaussian variables, where we use
complex variables whenever the underlying Banach spaces are complex.
Finally, let us fix the orthonormal systems $G_n:=(g_k)_{k=1}^n$ and
$U_n:=(e_k)_1^n$, where $(e_k)$ is the unit vector basis of $\lzn$. \hz

The starting point for our investigations is the following well-known
observation.

\begin{ttheorem} Let $T\in\comxy$ and let $\Phi = (\phi_k)_{k\in\nz}$ be an
orthonormal system.
\begin{enumerate}
\item[i)] If T is of cotype 2, then one has for all finite sequences $(x_k)
	  \subset X$
	  \[ \Tfmm \pl \le \pl c_2(T) \pl \gmm \pl. \]
\item[ii)] If T is of type 2, then one has for all finite sequences $(x_k)
	  \subset X$
	  \[ \Tgmm \pl \le \pl t_2(T) \pl \fmm \pl .\]
\end{enumerate}
\end{ttheorem}\hz

Let us recall that an operator $T\in\comxy$ is of cotype 2 or type 2 if there
are constants $c>0$ or $t>0$ such that for all finite sequences
$(x_k)\subset X$
\[ \kla \summ_k \| Tx_k\|^2 \mer^{\frac{1}{2}} \kl c \p \gmm
   \pll \mbox{or} \pll
   \Tgmm \kl t \p \kla \summ_k \| x_k\|^2 \mer^{\frac{1}{2}} \p . \]
As usual $c_2(T) := \inf c$ and $t_2(T) := \inf t$.
Restricting the above inequalities to $n$ vectors $(x_k)_1^n$ we obtain
$c_2^n (T)$ and $t_2^n (T)$ which can be defined for all $T\in\comxy$.\hz \\
For a further discussion of Theorem 1 as well as for the first
occurrence we refer to \cite{LE-TA}(Theorems 9.24, 9.25). Clearly, Theorem 1
expresses
the extreme position of gaussian variables among arbitrary orthonormal
systems.
An easy approximation argument yields a converse of Theorem 1 in the case of
complete orthonormal systems.

\begin{ttheorem} Let $\on_{k \in \nz}$ be a \underline{complete}
orthonormal system in $L_2(0,1)$ and let $T\in\comxy$.
If there is a $c >0$ such that for all finite sequences
$(x_k) \subset X$
	  \[ \Tfmm \pl \le \pl c \pl \gmm \p , \]
then T is of cotype 2.
If there is a $c >0$ such that for all finite sequences
$(x_k) \subset X$
	  \[ \Tgmm \pl \le \pl c \pl \fmm \p , \]
then T is of type 2.
\end{ttheorem}

$Proof:$ For example, one can use \p
  $\sum_{l=1}^N \| y_l\|^2 \p = \p \int_0^1 \noo \sum_{k=1}^{\infty} \kla
   \sum_{l=1}^N (f_l,\phi_k) y_l \mer \phi_k (t) \rrm^2 dt$ \p
for \linebreak
$f_l (t) := \sqrt{l(l+1)} \chi_{(\frac{1}{l+1},\frac{1}{l})}(t)$
and the fact that $g_l^{'} = \sum_{k=1}^{\infty} (f_l,\phi_k) g_k$ are again
standard gaussian variables which are independent. \hfill $\Box$ \hz

In this paper we mainly discuss the following two problems.

\begin{enumerate}
\item [(P1)] Characterize those (not necessarily complete) orthonormal systems
      $\Phi =(\phi_k)_{k\in\nz}$ such that the conclusions of Theorem 2
      hold true.
\item [(P2)] Find a local version of Theorem 2 in the sense that we consider
      systems $\fis$ of a given length and ask for the usual cotype and type
      constants restricted to $n$ vectors.
\end{enumerate}

To give a systematic treatment let $\fiss \subset L_2(\Om,p)$,
(I is a countable index set) be an orthonormal system. An operator
$T\in\comxy$ is said to be of $\Phi$-cotype and  $\Phi$-type, respectively,
if there is a constant $c\ge 0$ such that

\[\left\| \sum\limits_{k\in I} \phi_k T x_k \right\|_{L_2(Y)}
  \p\le \pl c \p \left\| \sum\limits_{k\in I} g_k x_k \right\|_{L_2(X)}
  \pl\mbox{and}\pll
  \left\| \sum\limits_{k\in I} g_k T x_k \right\|_{L_2(Y)}\p\le \pl c \p
  \left\| \sum\limits_{k\in I} \phi_k x_k \right\|_{L_2(X)} \pl ,\]

respectively. The best possible constants will be denoted by $\fco T$ and
$\fty T$. Furthermore, $T\in\comxy$ is said to be of the modified $\Phi$-type,
if there is a constant $c\ge 0$ such that for all $f\in L_2(X)$

\[ \noo \summ_{k\in I} g_k \kla\int Tf \bar{\phi}_k dp \mer \rrm_{L_2(Y)} \kl
   c \p \| f\|_{L_2(X)} \p .\]

The best possible constant is denoted by $\ftyh T$. In the same way one could
define the remaining case $\hat{c}_{\Phi}(T)$. Since this case follows by
duality from the modified $\bar{\Phi}$-type \linebreak
($\bar{\Phi}:=(\bar{\phi_k})$
is the conjugate system) we will omit this case. To consider the above
quantities we use ideal norms as one of the main tools. Given an
orthonormal system $\fis$ we first introduce for $u\in\com X)$
\[ \Phi (u) := \noo\summ_{k=1}^n \phi_k ue_k \rrm_{L_2(X)}. \]

The dual norm $\Phi^*$ on ${\cal L}(X,\lzn )$ is given by trace duality, that
is

\[ \Phi^* (v) \p := \p \sup \left\{ | {\rm tr} (vu) | \p \mitt \p \Phi (u) =1
   \right\} \p .\]

Note that in general the definition of $\Phi (u)$ and $\Phi^*(v)$ depends on
the special choice of the orthonormal system $\{ e_1 ,...,e_n\}\subset\lzn$
(we can consider $\Phi$ as a norm on $X^n$ via $\Phi ((x_1,...,x_n )) := \noo
\summ \phi_k x_k \rrm_
{L_2(X)}$ and $\Phi^*$ as the dual norm on $[X^n,\Phi ]^*$). There are two
standard procedures to generate ideal norms starting from $\Phi$.
For $T\in \comxy$ we define

\[ \pif (T) := \sup \left\{ \Phi (Tu) \pl | \pl
			    \| u :\lzn \rightarrow X\| \le 1 \right\}, \]
and, if $T$ is a finite rank operator,

\[ \nu_{\Phi} (T) := \inf \left\{ \summ_{j=1}^N \Phi (u_j) \|v_j\| \pl \mitt
			    \pl T =\summ_{j=1}^N u_j v_j , \pl u_j \in\com Y),
			    \pl v_j \in{\cal L}(X,\lzn ) \right\}. \]

It is an easy exercise to check that $\pif$ is an ideal norm on the class of
all bounded operators and that $\nu_{\Phi}$ is an ideal norm on the class
of finite rank operators. That is, we have for
$\al\in\{\pi_{\Phi} ,\nu_{\Phi}\}$  the norm properties  and the relations

\[ \| BTA\| \le \al (BTA) \le \| B\| \al (T) \| A\| \]

and $\al (a\otimes y) = \| a\|_{X^*} \| y\|_Y$ for $a\otimes y \in \comxy$.
To shorten the notation we will write in the sequel $\al (X)$ instead of
$\al (I_X)$. The connection to the approximation theory is given by the
geometric
interpretation of the inequalities $\pif (\lin )\ge \delta^{'} \sqrt{n}$
and $\pif (\izin )\ge \delta^{''} \sqrt{n}$. They correspond to the conditions
$(G_n^{'})$ and $(G_n^{''})$, respectively, in Theorem \ref{geometry} below.
To compare the usual type and cotype with the $\Phi$-type and cotype
we will compare the $\pi_{\Phi}$-norm with the $\pi_2^n$-norm directly.\hz \\
Let us remember that an operator $T\in\comxy$ is absolutely q-summing
($1\le q<\infty$) provided there is a constant $c>0$ such that for all finite
sequences $(x_k)\subset X$ one has

\[ \kla \summ_k \| Tx_k \|^q \mer^{\frac{1}{q}} \kl c \p
   \sup \left\{ \kla \summ_k |\langle x_k,a\rangle |^q \mer^{\frac{1}{q}} \p
		\mitt \p a\in B_{X^*} \p \right\} .\]

The best possible constant is denoted by $\pi_q(T)$. Considering the above
inequality for $n$ vectors $(x_k)_1^n$ only we get $\pi_q^n(T)$ which is
again defined for all $T\in\comxy$.\hz \\
The concept of the $\pi_{\Phi}$-norms connects in a natural way the usual
$\ell$-norm with the $\pi_2^n$-norm. Namely, if $G_n = (g_k)_1^n$ and
$U_n = (e_k)_1^n$ are defined as below then we recover the $\ell$-norm
and the $\pi_2^n$-norm for $u\in\com X)$ by

\[ \ell (u) \p = \p \pi_{G_n} (u)\p = \p\noo \summ_1^n g_k ue_k \rrm_{L_2(X)}
   \pla \mbox{and} \pla \pi_2^n (u) \p =\p \pi_{U_n} (u) \p .\]

These are the extreme situations since in any case

\[ \ell (u) \kl \pif (u) \kl \pi_2 (u) \kl \sqrt{2} \pi_2^n (u)
   \pla \mbox{whenever} \pla u\in\com X) \pl \]

(see Remark \ref{ell_pif} , Lemma \ref{phi_p2} and \cite{TOM} for the latter
inequality). The following example has served us as a pro type for the whole
investigation and also as a motivation for the introduction of $\pi_{\Phi}$
and $\nu_{\Phi}$ norms
in connection with the problems concerning
type and cotype.

\begin{example} Let $\Phi = E_n := \kla e^{ikt} \mer_{k=1}^n \subset
L_2 (\Pi)$ be the trigonometric system. Then

\[ c_2^n(T) \kl c \p c_{E_n}(T) \hspace{.7cm} \mbox{and} \hspace{.7cm}
   t_2^n(T) \kl c \p t_{E_n}(T) \hspace{.7cm} \mbox{for all} \hspace{.7cm}
   T\in\comxy , \]

whereas $c>0$ is an absolute constant independent from $n$.
\end{example} \hz

$Proof$: Using the Marcinkiewicz--Zygmund--inequality (see
\cite{ZYG}(II,p.30),
\cite{PIE-WEN})

\[ \kla \frac{1}{n} \summ_{l=1}^n \noo \summ_{k=1}^n e^{2\pi i \frac{kl}{n}}
   x_k \rrm^2_X \mer^{\frac{1}{2}} \kl c \p\kla\int\limits_0^{2\pi} \noo
   \summ_{k=1}^n  e^{ikt} x_k \rrm^2_X \frac{dt}{2\pi} \mer^{\frac{1}{2}} \]

implies, by a simple rotation argument, for all $u\in\com X)$

\begin{eqnarray*}
\pi_2^n (u)
& = & \sup_{\| w : \lzn \rightarrow \lzn \| =1}
       \left\{ \kla \summ_{k=1}^n \noo uwe_k \rrm^2\mer^{\frac{1}{2}} \right\}
    = \p \sup_{\| w : \lzn \rightarrow \lzn \| =1}
       \left\{ \kla \summ_{l=1}^n \noo \summ_{k=1}^n
	       \frac{e^{2\pi i \frac{kl}{n}}}{\sqrt{n}} \p uwe_k \rrm^2
	       \mer^{\frac{1}{2}} \right\} \\
&\le&  c \p \sup_{\| w : \lzn \rightarrow \lzn \| =1}
       \left\{ \kla \int\limits_0^{2\pi} \noo \summ_{k=1}^n
	       e^{ikt} \p uwe_k \rrm^2 \frac{dt}{2\pi}
	       \mer^{\frac{1}{2}} \right\}
    = \p c \p \pi_{E_n} (u) \p .
\end{eqnarray*}

Consequently,

\begin{enumerate}
\item[$(\ast )$] \hspace{3.5cm} $\pi_2^n (u) \kl c \p\pi_{E_n} (u)$ .
\end{enumerate}

Hence
$\pi_2^n (Tu) \kl c \p \pi_{E_n} (Tu) \kl c \p c_{E_n} (T) \p \ell (u)$ and
$c_2^n (T) \kl c \p c_{E_n} (T)$ for all $T\in\comxy$.
Let us turn to the type situation. To deduce the type equivalence
from $(\ast )$ we cannot use the Riesz-projections
since this would require UMD-properties (for example) for $X$
(see e.g. \cite{PIE-WEN}). Instead of this we use the de la
Vall\'{e}e Poussin kernel and find a sequence $\lambda_1 ,....,\lambda_n
\ge 0$ with $\lambda_i =1$ for $n_0 \le i\le n_1$, whereas $n_1 -n_0 \ge
\frac{n}{3}$, such that
\[ \noo \summ_{k=1}^n \lambda_k \hat{f} (k) e^{ikt} \rrm_{L_2(X)} \le c_1
   \| f\|_{L_2(X)} \]

for all $f\in L_2(X)$ and some absolute constant $c_1>0$ independent from $n$
($\hat{f}(k) = \int_0^{2\pi} f(t) e^{-ikt} dt/2\pi$).
Now let $1\kl r \kl \frac{n}{3}$ and $u\in{\cal L}(\ell_2^r ,X)$ be such that
$\hat{f}(k+N_0) = ue_k$ for $k=1,...,r$.
Defining $J:\ell_2^r \rightarrow \ell_2^n$ by $Je_k := e_{k+n_0-1}$ and
$\tilde{u}:\lzn \rightarrow X$ by $\tilde{u} := \summ_{k=1}^n \lambda_k e_k
\otimes \hat{f} (-n_0+1+N_0+k)$ we obtain $u=\tilde{u}J$ and
\[ \nu_{E_n} (u) \le E_n(\tilde{u}) \le c_1 \noo f e^{-i(-n_0+1+N_0)t} \rrm_2
   = c_1 \| f \|_2. \]
Using a simple blocking argument and the definition of $\nu_{E_n}$ this means
for $w\in \com X)$

\begin{enumerate}
\item[($\ast\ast$)] \hspace{1.7cm} $\nu_{E_n} (w) \le c_2 \inf\left\{ \|
	  f\|_{L_2(X)} \mitt \hat{f} (k) = we_k\p ,\p k=1,...,n \right\}$ \p
	  .
\end{enumerate}
This allows us to consider the bigger norm $\hat{t}_{E_n}$.
Indeed, we get
\[ \noo \summ_{k=1}^n T \hat{f}(k) g_k \rrm_{L_2(Y)}
   \kl t_{E_n}(T) \p \nu_{E_n} \kla \summ_{k=1}^n e_k\otimes\hat{f}(k)\mer
   \kl c_2 \p t_{E_n}(T) \p \| f\|_{L_2(X)} \]

such that $\hat{t}_{E_n} (T) \le c_2 t_{E_n} (T)$. Now we are in the position
to use a duality argument.
The Marcinkiewicz--Zygmund--inequality gives for ${g'}_l :=
\frac{1}{\sqrt{n}} \sum_{k=1}^n e^{2\pi i \frac{kl}{n}} \bar{g}_k$ and
$h\in L_2(Y^*)$

\begin{eqnarray*}
\kla \summ_{l=1}^n \noo \int T^* h\bar{g'}_l dp\rrm^2_{X^*}\mer^{\frac{1}{2}}
& = & \kla \summ_{l=1}^n \noo \summ_{k=1}^n \frac{e^{-2\pi i\frac{kl}{n}}}
      {\sqrt{n}} \int T^* h g_k dp \rrm^2_{X^*} \mer^{\frac{1}{2}} \\
&\le& c \p \kla \int\limits_0^{2\pi} \noo \summ_{k=1}^n  e^{-ikt}
      \int T^* h g_k dp \rrm^2_{X^*}\frac{dt}{2\pi}\mer^{\frac{1}{2}}\p .
\end{eqnarray*}

Using the convenient duality properties of the modified type we can continue
with

\[ \kla \int\limits_0^{2\pi} \noo \summ_{k=1}^n  e^{-ikt}
   \kla \int T^* h g_k dp \mer \rrm^2_{X^*}\frac{dt}{2\pi}\mer^{\frac{1}{2}}
   \kl \hat{t}_{E_n} (T) \p \| h\|_{L_2(Y^*)} \p . \]

Consequently,

\[ \kla \summ_{l=1}^n \noo \int T^* h \bar{g'}_l dp \rrm^2_{X^*} \mer^
   {\frac{1}{2}} \kl cc_2 \p t_{E_n} (T) \| h\|_{L_2(Y^*)} \p . \]

Using again duality we arrive at $t_2^n (T) \kl c c_2 t_{E_n} (T)$.
\hfill $\Box$ \hz

The paper is organized in the following way. First we consider the problem
(P2) mentioned above and derive as a simple consequence the answer of (P1).
Concerning the problem (P2) our main theorem states that it is sufficient
to test cotype and type conditions on rather extreme operators. More
precisely we prove

\begin{ttheorem} \label{phi-type} Let $\fis$ be an
orthonormal system and let $\delta >0$. For some absolute constant $c>0$ not
depending on $\delta$,$n$, and $\Phi$, the following holds true.
\begin{enumerate}
\item If $\fco \lin \ge \delta \sqrt{\frac{n}{\log (n+1)}}$
      or $\ftyh \len \ge \delta \sqrt{n}$ then
      $\delta^3 \pz w) \le c \pif (w)$ for all $w\in \com X)$.
      Consequently, for all $T\in {\cal L}(X,Y)$
      \[ \frac{\delta^3}{c} c_2^n (T) \pl \le \pl \fco T
	 \pl \le \pl \sqrt{2} c_2^n (T) \pll \mbox{and} \pll
	 \frac{\delta^3}{c} t_2^n (T) \pl \le \pl \ftyh T \pl \le
	 \pl \sqrt{2} t_2^n (T).\]
\item If $\fty \len \ge \delta \sqrt{n}$ then
      $\delta^3 \vf (w) \le c \pz w^*)$ for all $w\in \com X)$.
      Consequently, for all $T\in {\cal L}(X,Y)$
      \[ \frac{\delta^3}{c} t_2^n (T) \pl \le \pl \fty T \pl \le  \pl
	 \ftyh T \pl \le \pl \sqrt{2} t_2^n (T). \]
\end{enumerate}
\end{ttheorem} \hz

It turns out that the conditions ($\ast$) and ($\ast\ast$) from Example 3
are necessary in general. Note, that $\pz u^*) \kl \Phi^*(u^*) \kl \| f\|_2$
if $\hat{f} (k) \p = \p  ue_k$ ($k=1,...,n$) (see section \ref{ons}) such that
Theorem \ref{phi-type}(2) implies ($\ast\ast$) with $c_2=\frac{c}{\delta^3}$.
In section \ref{abs} we establish an abstract version of Theorem
\ref{phi-type} in the terms of operator ideals whereas the connection to the
notion of cotype and type is given in section \ref{ons}.
The proof of Theorem \ref{phi-type}(1) consists of several steps
formulated in the next theorem which is verified in section \ref{ons}.
In the first step $(G_n) \nach (G_n^{'})$ one we get rid of the logarithmical
factor but lose the orthogonality. In the second one $(G_n^{'}) \nach
(G_n^{''})$, which is an essential part of the proof of Theorem
\ref{phi-type},
we come back to an orthonormal system. The last condition in the abstract
corresponds to $(G_n^{''})$ via an observation of Bourgain.


\begin{ttheorem} \label{geometry} Let $\fis \subset L_2(\Om ,p)$ be an
orthonormal system and let $\delta ,\delta^{'},\delta^{''}>0$. Let us define
the following conditions.
\begin{enumerate}
\item[$(G_n)$] $\fco \lin \ge \delta \sqrt{\frac{n}{\log (n+1)}}$\p .
\item[$(G^{'}_n)$] There exists $n$ functions $h_j \in
     {\rm span}\{\phi_k \p |\p k=1,...,n\}$ with $\noo h_j \rrm_2 \kll 1$
     and
     \[ \kla \int\limits_{\Om} \sup_{j=1,...,n} \bet h_j \rag^2 \pl dp
	\mer^{\frac{1}{2}} \gl \delta^{'} \pl \sqrt{n} \pl .\]
\item[$(G^{''}_n)$] There exists an orthonormal system $\Psi =(\psi_j)_1^n
     \subset {\rm span} \{\phi_k \p |\p k=1,...,n\}$ with
     \[ \kla \int\limits_{\Om} \sup_{j=1,...,n} \bet \psi_j \rag^2 \pl dp
	\mer^{\frac{1}{2}} \gl \delta^{''} \pl \sqrt{n} \pl .\]
\end{enumerate}
Then
$(G_n) \rightarrow (G^{'}_n) \rightarrow (G^{''}_n) \rightarrow (G_n)$
with $\delta^{'} = \frac{\delta}{20c_0 \sqrt{1+\log
(\frac{c_0}{\delta} +1)}}$,
$\delta^{''} = \frac{1}{c_1} \delta^{'3}$, and $\delta = \frac{1}{c_2}
\delta^{''}$ for numerical constants $c_0, c_1, c_2 >0$.
\end{ttheorem}

In section \ref{few} we consider orthonormal systems defined on measure spaces
with
few atoms and no continuous part. We prove the following Theorem \ref{sumfew}
which uses
the local theory of Banach spaces to clarify the relations between cotype and
type conditions.


\begin{ttheorem} \label{sumfew} Let $1\kll n \kll N$. Then one has the
following.
\begin{enumerate}
\item [(1)] $c_2^n(T) \kl 12 \pl \sqrt{\frac{N}{n}}\pl
      \fco T$ \p for all $T\in\comxy$ and all $\fis \subset \ell_2^N$ \p.
\item [(2)] For all $2<q<\infty$ there is an orthonormal system $\fis \subset
      \ell_2^N$ such that
      \[ \fco \lin \kl c_q \p \max \left\{\sqrt{q}\, n^{1/q-1/2}, \kla
	 \frac{n}{N} \mer^{1/q} \right\} \p c_2^n(\lin )\p. \]
\item [(3)] As long as $0<\vare <1$ and $1\kll n \kll (1-\vare)N$ there is
      an orthonormal system $\fis \subset \ell_2^N$ with
      \[ \fty {\ell_1} \kl c_0\pl \sqrt{\frac{1}{\vare}\p\log \kla
      1+\frac{1}{\vare}\mer}\pl.\]
\end{enumerate}
$c_0>0$ is an absolute constant whereas $c_q>0$ depends on $q$ only.
\end{ttheorem}

Assertion (1) shows that as long as $N$ is proportional to $n$ the
corresponding cotype constants are equivalent, whereas in (3) "pathological"
orthonormal systems are found for the notion of type. The theory of
$\Lambda_p$-sets is involved for the construction of the orthonormal systems
in the second assertion. Choosing $n \sim N^{\delta}$ we obtain systems which
fail the first conclusion of Theorem \ref{sumfew}.


\section{Abstract theory}\label{abs}

Throughout this section we will say that a norm $\al$ on
$\com \cdot )$ (which means the collection of all $\com X)$, where $n$ is
fixed and $X$ is an arbitrary Banach space) is an ideal norm if

\[ \| T u A \| \kl \al (T u A) \kl \| T\| \p \al (u) \p \| A\| \hspace{.6cm}
   \mbox{and}\hspace{.6cm}\al (a\otimes x)\pl =\p\| a\|_{\lzn} \p \| x\|_X \]

for all $T\in\comxy$\p ,\p $u\in\com X)$\p ,\p $A\in\com\lzn )$\p ,\p
$a\in\lzn$ \p and $x\in X$. Similarly, $\beta$ is an ideal norm on
${\cal L}(\cdot ,\lzn )$ if

\[ \| A v T \| \kl \beta (A v T) \kl \| A\|\p\beta (v)\p \| T\| \hspace{.6cm}
   \mbox{and} \hspace{.6cm} \beta (b\otimes y) \p = \p \| b\|_{Y^*} \p
   \| y\|_{\lzn} \]

for all $A\in\com\lzn )$\p ,\p $v\in{\cal L}(Y,\lzn )$\p ,\pl $T\in\comxy$\p ,
\pl $b\in Y^*$ and $y\in\lzn$. The adjoint ideal norms $\al^*$ on
${\cal L}(X,\lzn )$ and $\beta^*$ on $\com X)$ are given by

\[ \al^*(v) \p = \p \sup_{\al (u:\lzn\nach X)\le 1} |{\rm tr} (vu)|
   \plll \mbox{and}\plll
   \beta^*(u) \p = \p \sup_{\beta (v:X\nach\lzn) \le 1} |{\rm tr} (vu)| \p .\]

Furthermore, we define $\epz T) := \inf c$ such that

\[ \kla \summ_{i=1}^n \| Tx_i\|^2 \mer^{\frac{1}{2}} \le c \sup \left\{
   \summ_{i=1}^n |\langle x_i ,a\rangle|^2 \mitt a\in B_{X^*} \right\} \]

for all $x_1 ,...,x_n \in X$ with $\| Tx_1\| =...=\| Tx_n\|$.
The following lemma is the key for what follows. The proof is similar to the
proof of \cite{BKT}(Theorem 3.1). We thank Th. K\H{u}hn for his hints to
improve the constant appearing in Lemma \ref{epz} . \hz

\begin{lemma} \label{epz} Let $T \in {\cal L}(X,Y)$ and $n \in \nz$. Then
\[ \pzn T) \le \sqrt{6} \pl\epz T). \]
\end{lemma}

$Proof:$ Assume $x_1 ,...,x_n \in X$ with $\summ_1^n \| Tx_i \|^2 =1$
whereas $\| Tx_i \| >0$ for all i. Setting

\[ \sigma_m := \left\{ i\in \{ 1,...,n\} \pl | \pl 2^{-m} < \| Tx_i \|^2
   \le 2^{1-m} \right\} \]

we obtain $\summ_{m=1}^{\infty} |\sigma_m | = n$. Let $m_0 \in \nz$ such
that $2^{m_0-1} < 3n \le 2^{m_0}$. Then we get

\[\summ_{m=m_0 +1}^{\infty} \summ_{i\in \sigma_m} \| Tx_i \|^2
  \le \summ_{m=m_0 +1}^{\infty} |\sigma_m | 2^{1-m}
  \le 2^{-m_0} \summ_{m=m_0 +1}^{\infty} |\sigma_m | \le 2^{-m_0} n
  \le \frac{1}{3}. \]

Now we define

\[ I:=\left\{ (i,j) \mitt j=1,...,2^{m_0-m} \pll \mbox{if} \pll i\in\sigma_m
	      \pll ; \pll m = 1,...,m_0 \right\} \]

and obtain

\[ |I| = \summ_{m=1}^{m_0} |\sigma_m | 2^{m_0 -m} \le 2^{m_0}
   \summ_{m=1}^{\infty} \summ_{i\in\sigma_m} \| Tx_i \|^2 \le 2^{m_0} < 6n \]

as well as

\begin{eqnarray*}
|I|
&  =  & \summ_{m=1}^{m_0} |\sigma_m | 2^{m_0 -m} \pl\ge\pl 2^{m_0 -1}
	\summ_{m=1}^{m_0} \summ_{i\in\sigma_m} \| Tx_i \|^2 \\
& \ge & 2^{m_0 -1} \kla 1- \summ_{m=m_0 +1}^{\infty} \summ_{i\in\sigma_m}
	\| Tx_i \|^2 \mer \pl\ge\pl 2^{m_0 -1} \kla 1- \frac{1}{3} \mer
	\pl\ge\pl n .
\end{eqnarray*}

Defining $y_{ij} := \frac{x_i}{\| Tx_i\|}$ for $(i,j)\in I$ and choosing a
subset $J\subseteq I$ with $|J| =n$ we deduce

\begin{eqnarray*}
n \pl = \pl |J|
& \le & \epz T)^2 \sup_{a\in B_{X^*}} \kla \summ_J |<y_{ij},a>|^2 \mer \\
& \le & \epz T)^2 \sup_{a\in B_{X^*}} \kla \summ_I |<y_{ij},a>|^2 \mer \\
& \le & \epz T)^2 \sup_{a\in B_{X^*}} \kla \summ_{m=1}^{m_0}
	\summ_{i\in\sigma_m} \summ_{j=1}^{2^{m_0-m}}
	|<\frac{x_i}{\| Tx_i\|},a >|^2 \mer \\
& \le & 2^{m_0} \epz T)^2 \sup_{a\in B_{X^*}} \kla \summ_{m=1}^{m_0}
	\summ_{i\in\sigma_m} |<x_i,a >|^2 \frac{2^{-m}}{\| Tx_i \|^2} \mer \\
& \le & 2^{m_0} \epz T)^2 \sup_{a\in B_{X^*}} \kla \summ_{m=1}^{m_0}
	\summ_{i\in\sigma_m} |<x_i,a >|^2 \mer \\
& \le & 2^{m_0} \epz T)^2 \sup_{a\in B_{X^*}} \kla \summ_{i=1}^n |<x_i,a >|^2
	\mer .
\end{eqnarray*}

Consequently,
\for \pzn T)\le\sqrt{\frac{2^{m_0}}{n}}\epz T)\le\sqrt{6}\epz T) \p.
\\[-1.5cm] \mel \hfill $\Box$\hz

\begin{lemma} \label{max_min} Let $\al$ and $\beta$ be ideal norms on
$\com \cdot)$ and ${\cal L}(\cdot ,\lzn )$, respectively. Then for all
$u\in\com X)$ and $v\in{\cal L}(X,\lzn )$
\begin{enumerate}
\item [(1)] $\al (\izin ) \le\al (u)$ \pll whenever \pll $\| ue_i\| =1$ for
      \pl $i=1,...,n$,
\item [(2)] $\beta (\iezn ) \le\beta (v)$ \p\pll whenever \pll $\| v^*e_i\|
=1$
      for \pl $i=1,...,n$.
\end{enumerate}
\end{lemma} \hz

$Proof:$ (1) Choosing $a_1,...,a_n \in B_{X'}$ with $\langle ue_i,a_i \rangle
=1$ and setting $w:= \summ_{i=1}^n a_i \otimes e_i \in {\cal L} (X,\lin )$ we
obtain

\[ n= {\rm tr}(\iizn wu) \le \al^* (\iizn ) \| w\| \al (u) \le
   \al^* (\iizn ) \al (u). \]

Using $\al (\izin ) \al^* (\iizn ) =n$ from \cite{PIEOP} (9.1.8) we arrive at
our assertion.
(2) For $\vare >0$ we choose $x_1,...,x_n \in B_X$ with
$\langle x_i,v^*e_i \rangle \ge 1-\vare$ and set
$w:= \summ e_i \otimes x_i \in {\cal L}(\len ,X)$. Hence

\[ (1-\vare )n \le {\rm tr}(\izen vw) \le \beta^* (\izen ) \beta (v) \| w\| \]

and $n\le \beta^* (\izen ) \beta (v)$ such that we finish as in (1).
\hfill $\Box$ \hz

\begin{lemma} \label{pre_theorem} Let $\al$ and $\beta$ be ideal norms on
$\com \cdot)$ and ${\cal L}(\cdot ,\lzn )$, respectively. Then
\begin{enumerate}
\item [(1)] $\al (\izin ) \pl \epz u) \le \sqrt{n} \pl \al (u)$ \pll for all
      \pll $u\in \com X)$,
\item [(2)] $\beta (\iezn ) \pl \epz v^* ) \le \sqrt{n} \pl \beta (v)$ \pll
      for all \pll $v\in {\cal L}(X,\lzn )$.
\end{enumerate}
\end{lemma} \hz

$Proof:$ (1) Let $w\in \com \lzn )$ be such that $\| uwe_i\| =1$ for
$i=1,...,n$. Then $\al (\izin ) \le \al (uw) \le \al (u) \| w\|$ and

\[ \al (\izin ) \sup \left\{ \frac{\sqrt{n}}{\| w\|} \pl \mitt \pl \| uwe_i
   \| =1 \hspace{.2cm} \mbox{for} \hspace{.2cm} i=1,...,n \right\} \le
   \sqrt{n} \pl \al(u) .\]

(2) Let $w\in\com \lzn)$ be such that $\| v^* w^* e_i\| =1$ for $i=1,...,n$.
Then $\beta (\iezn ) \le \beta (wv) \le \| w\| \beta (v)$ and

\for \beta (\iezn ) \sup \left\{ \frac{\sqrt{n}}{\| w^* \|} \pl \mitt \pl \|
   v^* w^* e_i \| =1 \hspace{.2cm} \mbox{for} \hspace{.2cm} i=1,...,n \right\}
   \le \sqrt{n} \pl \beta (v) \p . \\[-1.5cm] \mel \hfill $\Box$ \hz

Combining Lemmata \ref{epz}, \ref{pre_theorem}, and the fact that
$\pz T) \le \sqrt{2} \pzn T)$ whenever ${\rm rank} (T) \le n$ (see \cite{TOM})
we get \hz

\begin{theorem} \label{main_abs} Let $\al$ and $\beta$ be ideal norms on
$\com \cdot)$ and ${\cal L}(\cdot ,\lzn )$, respectively. Then
\begin{enumerate}
\item [(1)] $\al (\izin ) \pz u) \le \sqrt{12} \pl \sqrt{n} \pl
      \al (u)$ \pll for all \pll $u\in\com X)$,
\item [(2)] $\beta (\iezn ) \pz v^*) \le \sqrt{12} \pl \sqrt{n}
      \pl \beta (v)$ \pll for all \pll $v\in{\cal L}(X,\lzn )$.
\end{enumerate}
\end{theorem} \hz

Let us recall that the $n$-th approximation number \cite{PIE} of an operator
$T\in\comxy$ is defined by

\[ a_n (T) \p := \p \inf\left\{ \| T-L \| \p \mitt \p L\in\comxy \p , \p
				{\rm rank} (L) < n \right\}\p . \]

To bring the above theorem in a form we need we will use

\begin{lemma} \label{lemma_delta^3} Let $\beta$ be a norm on
${\cal L} (\len ,\lzn )$ such that for all $v\in{\cal L} (\len ,\lzn )$ and
all orthogonal matrices $w\in\com \lzn )$ one has
$\be (wv) = \be (v) \le \nu (v)$. If $\sup\left\{\beta (v) \pl |\pl
\| v:\len \rightarrow \lzn \| \le 1 \right\} \ge \delta \sqrt{n}$ for some
$\delta >0$, then $\beta (\iezn )\ge \frac{\delta^3}{c} \sqrt{n}$ where $c >0$
is an absolute constant.
\end{lemma}

$Proof:$ Using Grothendieck's inequality (see \cite{PIS}(Theorem 5.10)) our
assumption ensures the existence of some $v\in{\cal L}(\len ,\lzn )$ with
$\pz v) \le K_G$ and $\beta (v) \ge \delta \sqrt{n}$. Trace duality gives
some $u\in \com \len )$ with $\pz u) \ge \frac{\delta}{K_G} \sqrt{n}$ and
$\beta^* (u) \le 1$ (note that $\be\le\nu$ implies $\|\cdot\| \le \be^*$).
Exploiting \cite{PIE} (2.7.4) we deduce for $\theta >0$

\[ \frac{\delta}{K_G} \sqrt{n} \le \pz u) \le c \summ_{k=1}^n
   \frac{a_k(u)}{\sqrt{k}} \le 2c \kla \sqrt{[\theta n]} + \sqrt{n}
   a_{[\theta n]+1} (u) \mer \]

and $\frac{\delta}{2cK_G} \le \sqrt{\theta} + a_{[\theta n]+1}(u)$. Setting
$\theta := \kla \frac{\delta}{4cK_G}\mer^2$ we obtain (using \cite{PIE}
(2.11.6,2.11.8)) for some orthogonal $w\in\com\lzn )$

\begin{eqnarray*}
\frac{\delta}{4cK_G}
&\le& a_{[\theta n]+1}(u)
 \le \sqrt{n} a_{[\theta n]+1}(\iezn u)
 \le \frac{\sqrt{n}}{[\theta n]+1} \summ_1^n a_k(\iezn u)
  =  \frac{\sqrt{n}}{[\theta n]+1} |{\rm tr} (w\iezn u)| \\
&\le& \frac{1}{\theta \sqrt{n}} \beta (w \iezn ) \beta^* (u)
  =  \frac{1}{\theta \sqrt{n}} \beta (\iezn ) \beta^* (u) \p .
\end{eqnarray*}

Hence $\frac{\delta^3}{(4cK_G)^3} \sqrt{n} \le \beta (\iezn )$ .
\hfill $\Box$ \hz

Now Theorem \ref{main_abs} and Lemma \ref{lemma_delta^3} imply \hz

\begin{cor} \label{cor_delta^3} Let $\al$ and $\beta$ be ideal norms on
$\com \cdot)$ and ${\cal L}(\cdot ,\lzn )$, respectively. Then
\begin{enumerate}
\item [(1)] $\sup \left\{ \frac{\al (w)}{\sqrt{n}} \pl \mitt \pl \| w:\lzn
      \rightarrow \lin \| =1 \right\}^3 \pz u) \pl \le \pl c \pl \al (u)$
      \pll for all \pll $u\in\com X)$,
\item [(2)] $\sup \left\{ \frac{\beta (w)}{\sqrt{n}} \pl \mitt \pl \| w:\len
      \rightarrow \lzn \| =1 \right\}^3 \pz v^*) \pl \le \pl c \pl \beta (v)$
      \pll for all \pll $v\in{\cal L}(X,\lzn )$,
\end{enumerate}
where $c>0$ is an absolute constant.
\end{cor} \hz

$Proof:$ (1) Setting $\beta (v):=\al (v^*)$ for $v\in {\cal L}(X,\lzn )$
and $\delta >0$ such that

\[ \delta \sqrt{n}
   = \sup \left\{ \al (w) : \| w:\lzn \rightarrow \lin \| =1 \right\}
   = \sup \left\{ \beta (v) : \| v:\len \rightarrow \lzn \| =1 \right\} \]

we obtain from Theorem \ref{main_abs} and Lemma \ref{lemma_delta^3}

\[ c \sqrt{n} \al (u) \ge \al (\izin ) \pz u) \ge \beta (\iezn ) \pz u)
   \ge \frac{\delta^3}{c'} \sqrt{n} \pz u). \]

Consequently, $\delta^3 \pz u) \le c c' \al (u)$.
(2) follows directly. \hfill $\Box$ \hz

In the following Corollary \ref{cor_delta^3} is made applicable to our
problems concerning type and cotype with respect to arbitrary orthonormal
systems. To do this we need the Weyl numbers and nuclear operators.
The $n$-th Weyl number \cite{PIE} of an operator $T\in\comxy$ is given by

\[ x_n(T) \p := \p \sup\left\{ a_n(Tu)\p\mitt\p u\in\com X)\p ,\p \| u\| = 1
			       \right\} \p .\]

An operator $T\in\comxy$ is nuclear \cite{PIEOP} provided that $T$ can be
written as

\[ T\p =\p \summ_1^{\infty} a_n\otimes y_n \]

with $a_n\in X^*$, $y_n\in Y$, and $\summ_1^{\infty}\| a_n\|\| y_n\|
<\infty$.
We set $\nu (T) := \inf \summ_1^{\infty}\| a_n\|\| y_n\|$ where the infimum is
taken over all possible representations.

\begin{lemma} \label{l-norm_to_op-norm} Let $\al$ be an ideal norm on
$\com \cdot )$ and let $u\in \com \lin )$ be such that
\[ \ell (u) \le 1 \hspace{1cm} \mbox{and} \hspace{1cm} \al (u) \ge \delta
   \sqrt{\frac{n}{\log (n+1)}} \]
for some $\delta >0$. Then there exists an operator $\tilde{u}\in \com \lin )$
with
\[ \| \tilde{u}\| = 1 \hspace{1cm} \mbox{and} \hspace{1cm} \al (\tilde{u})
   \ge \frac{\delta}{c_0 \sqrt{A}} \sqrt{n} \]
for all $A>1$  whenever $n\ge \kla \frac{c_0 \sqrt{A}}{\delta}
\mer^{\frac{2A}{A-1}}$. Moreover
\[ \al (\tilde{u}) \ge \frac{\delta} {20 \p c_0 \p \sqrt{1+\log
   (\frac{c_0}{\delta}+1)}} \sqrt{n} \]
for $n=1,2,...$. The constant $c_0>0$ is independent from $n,\delta ,A$ and
$\al$.
\end{lemma} \hz

$Proof:$ First we observe that
$a_r(u) \le c \frac{\ell (u)}{\sqrt{\log (r+1)}}$ for $u\in\com\lin )$,
which follows for example from the much deeper factorization
$u=BDA$ due to Talagrand used in the proof of Lemma \ref{Kaq}. This gives
the existence
of an orthogonal projection $P\in \com \lzn )$ with ${\rm rank} (P) > n-r$ and
$\| uP\| = a_r(u) \le c \frac{1}{\sqrt{\log (r+1)}}$. Hence via trace duality
we find an operator $v\in {\cal L}(\lin ,\lzn )$ with $\al^* (v)=1$ and

\begin{eqnarray*}
\delta \sqrt{\frac{n}{\log (n+1)}}
&\le& |{\rm tr}(vu)| \le |{\rm tr}(vuP)| + |{\rm tr}(vu(I-P))| \\
&\le& \nu (u) \| uP\| + 2 \summ_{k=1}^{r-1} x_k(v) a_k(u(I-P)) \\
&\le& \nu (u)  a_r(u) + 2 \summ_{k=1}^{r-1} x_k(v) a_k(u).
\end{eqnarray*}

Since Grothendieck's inequality \cite{PIS} implies

\[ x_k(v) \le k^{-1/2} \pz v) \le K_G k^{-1/2} \| v\|
   \le K_G k^{-1/2} \al^*(v) \le K_G k^{-1/2} \]

we can continue to

\begin{eqnarray*}
\frac{\delta}{c} \sqrt{\frac{n}{\log (n+1)}}
&\le& \frac{\nu (v)}{\sqrt{ \log (r+1)}} + 2 K_G \summ_{k=1}^{r-1}
      \frac{1}{\sqrt{k \log (k+1)}} \\
&\le& \frac{\nu (v)}{\sqrt{ \log (r+1)}} + 2 K_G c'
      \sqrt{\frac{r}{\log (r+1)}}.
\end{eqnarray*}

Hence, for $1\le r\le n$,

\[ \frac{\delta}{c} \sqrt{\frac{\log (r+1)}{\log (n+1)}} \sqrt{n} -
   2c' K_G \sqrt{r} \le \nu (v). \]

Now we pick for $A>1$ an $r\in \nz$ with $1\le r\le n$ and
$r\le (n+1)^{\frac{1}{A}} \le r+1$ and obtain

\[ \frac{\delta}{c\sqrt{A}} n^{\frac{1}{2}} - 4c' K_G
   n^{\frac{1}{2A}} \le \nu (v) .\]

Consequently, $n\ge\kla\frac{8 c c' K_G\sqrt{A}}{\delta}\mer^{\frac{2A}{A-1}}$
implies $\nu (v) \ge \frac{\delta}{2c} \frac{1}{\sqrt{A}} n^{\frac{1}{2}}$.
Finally, the desired operator $\tilde{u}\in \com \lin )$ is chosen such that
$\| \tilde{u}\| =1$ and $|{\rm tr}(v\tilde{u} )|=\nu (v)$. Setting
$c_0 :=\max (8cc'K_G,2c)$ we arrive at the first part of our assertion.
To prove the second assertion we put
$A_0 :=2+2\log \kla\frac{c_0}{\delta} +1\mer \ge 2$. It is clear that it
remains to consider the situation $\frac{c_0\sqrt{A_0}}{\delta}\ge1$ and
$n < \kla \frac{c_0\sqrt{A_0}}{\delta} \mer^{\frac{2A_0}{A_0 -1}}$. Here we
use $\frac{A_0}{A_0 -1}\le 1+\frac{2}{A_0}$ and
$\kla \frac{c_0\sqrt{A_0}}{\delta} \mer^{\frac{2}{A_0 }} \kl e^2$ to conclude
(for any $\tilde{u}$ with $\|\tilde{u}\| = 1$)
\for \frac{\delta}{e^2 c_0 \sqrt{A_0}} \sqrt{n}
   \kl \kla \frac{\delta}{c_0 \sqrt{A_0}}\mer^{1+\frac{2}{A_0}} \sqrt{n}
   \kl \kla \frac{\delta}{c_0 \sqrt{A_0}}\mer^{\frac{A_0}{A_0-1}} \sqrt{n}
   \kl 1 \kl \al (\tilde{u}) \\[-1.5cm] \mel \hfill $\Box$ \hz

\pagebreak

The main result of this section is \hz

\begin{theorem} \label{alpha_pz} Let $\al$ be an ideal norm on $\com \cdot )$
and let $\delta >0$.
\begin{enumerate}
\item [(1)] If there is an operator $u\in \com \lin )$ such that $\ell (u)
      \le 1$ and $\al (u) \ge \delta \sqrt{\frac{n}{\log (n+1)}}$ then
      \[ \delta^3 \pz w) \le c \al (w) \hspace{1cm}
      \mbox{for all} \hspace{1cm} w\in \com X). \]
\item [(2)] If there is an operator $u\in \com \len )$ such that $\al (u)
      \le 1$ and $\ell (u) \ge \delta \sqrt{n}$ then
      \[ \delta^3 \al (w) \le c \pz w^*)
      \hspace{1cm} \mbox{for all} \hspace{1cm} w\in \com X). \]
\end{enumerate}
$c>0$ is an absolute constant independent from $n,\delta$ and $\al$.
\end{theorem} \hz

$Proof:$ (1) Corollary \ref{cor_delta^3} and Lemma \ref{l-norm_to_op-norm}
imply (for some $c,c_0>0$)
$\kla \frac{\delta}{c_0 \sqrt{A}}\mer^3 \pz w) \le c \al (w)$ for all $w\in
\com X)$ whenever $n\ge \kla\frac{c_0\sqrt{A}}{\delta}\mer^{\frac{2A}{A-1}}$.
Setting $A=3/2$ we obtain
\[\pz w) \le \sqrt{n} \| w\| \le \sqrt{n} \al (w) \le \kla \frac{c_0 \sqrt{A}}
  {\delta}\mer^3 \al (w) \]
in the case $n< \kla \frac{c_0 \sqrt{A}}{\delta} \mer^{\frac{2A}{A-1}}$.
(2) Using trace duality we find
$v\in {\cal L} (\len ,\lzn )$ with $\ell^* (v) \le 1$ and $\al^* (v) \ge
\delta
\sqrt{n}$. Applying \cite{PS3} and \cite{TOM}(Theorem 12.7) we get
$\ell (v^*) \le K \sqrt{\log (n+1)}$ for some numerical constant $K>0$.
Lemma \ref{l-norm_to_op-norm} (applied to $\tilde{\al}(w) = \al^*(w^*)$)
produces an operator $\tilde{u} \in \com \lin )$ with $\|\tilde{u}\|\le 1$ and
$\al^* (\tilde{u}^*) \ge \frac{\delta}{Kc_0\sqrt{A}} \sqrt{n}$ whenever
$n\ge\kla \frac{Kc_0\sqrt{A}}{\delta}\mer^{\frac{2A}{A-1}}$.
Corollary \ref{cor_delta^3} (applied to $\be (w) = \al^* (w)$) shows
$\kla \frac{\delta}{Kc_0 \sqrt{A}} \mer^3 \pz w^*) \le c \al^*(w)$ for all
$w\in {\cal L}(X, \lzn )$. Now trace duality gives for $\tilde{w} \in \com X)$
the existence of some $w\in {\cal L}(X,\lzn )$ with $\al^*(w)=1$ and
\[ \al (\tilde{w}) = {\rm tr}(w\tilde{w}) \le \pz w^*) \pz \tilde{w}^*)
   \le \al^* (w) c \kla \frac{Kc_0\sqrt{A}}{\delta} \mer^3 \pz \tilde{w}^*)
   \le c \kla \frac{Kc_0\sqrt{A}}{\delta} \mer^3 \pz \tilde{w}^*). \]
In the case $n< \kla \frac{c_0 \sqrt{A}}{\delta} \mer^{\frac{2A}{A-1}}$ we
can continue as in (1) since $\al (\tilde{w}) \le \nu (\tilde{w}) \le
\sqrt{n} \pz \tilde{w}^*)$. \hfill $\Box$
\pagebreak


\section{Orthonormal systems in connection with type and cotype}\label{ons}
\setcounter{lemma}{0}

To handle the $\Phi$-type and cotype norms it is sometimes convenient to
introduce for orthonormal systems $\fis$ and $\psis$ and an operator
$T\in\comxy$ the quantity $\fipsih T := \inf c$ , such that

\[ \noo\summ_{k=1}^n \phi_k \kla \int Tf \bar{\psi}_k dp \mer \rrm_{L_2(Y)}
   \pl \le \pl c \pl \| f\|_{L_2(X)} \]

for all $f\in L_2(X)$ (see \cite{PIE-WEN}). It is clear that $\fipsih T =
\delta (T^*|\bar{\Psi},\bar{\Phi})$ whereas $\bar{\Phi} := (\bar{\phi_k})_1^n$
is the conjugate system of $\fis$. Using the same arguments as in
\cite{FI-TOM} (Lemma 9.2) (\cite{TOM} (Theorem 12.7)) it turns out that
\[ \fipsih T = \sup \left\{\Phi (Tu)\pl |\pl\bar{\Psi}^* (u^*) \le 1, u\in\com
   X) \right\},\]
where $\bar{\Psi}^*$ stands for $(\bar{\Psi})^*$.
In fact, for $u\in\com X)$ one obtains

\begin{eqnarray*}
\bar{\Psi}^* (u^*)  & = &  \sup\left\{ \left | \summ_1^n \langle ue_k ,
       a_k\rangle \right | \pl \mitt \pl \noo \summ_1^n \bar{\psi}_k a_k
       \rrm_2  \le 1 \right\}\\
& = &  \sup\left\{ \left | \int \left\langle f(\om ),\summ_1^n \bar{\psi}_k
       (\om )  a_k \right\rangle dp(\om ) \right | \pl \mitt \pl \noo
       \summ_1^n \bar{\psi}_k a_k \rrm_2 \le 1\p , \p
       \int f \bar{\psi_k} dp = ue_k \p \right\} \\
& = &  \inf \left\{ \noo f : {\rm span}\left\{\summ_1^n \bar{\psi}_k a_k
       \right\} \nach \kz \rrm\p
       \mitt \p f\in L_2(X) \p , \p  \int f \bar{\psi_k} dp = ue_k \p \right\}
\end{eqnarray*}

and

\[ \fipsih T \kl \sup \left\{\Phi (Tu)\pl |\pl\bar{\Psi}^* (u^*) \le 1,
   u\in\com X) \right\}.\]
For the reverse inequality we take $g\in L_2(Y^*)$ with $\| g\|_2 \kl 1+\vare$
and $\Phi (Tu) = \langle \summ_1^n \phi_k Tue_k,g\rangle$ such that
\begin{eqnarray*} \Phi (Tu)
& = & \int \left\langle\summ_{k=1}^n \psi_k ue_k ,\summ_{l=1}^n \bar{\psi}_l
      \int \phi_l T^*gdp \right\rangle dp \\
&\le& \noo\summ_{k=1}^n \psi_k ue_k : {\rm span} \left\{ \summ_1^n
      \bar{\psi}_k a_k \right\} \nach \kz \rrm \p \delta (T^*|\bar{\Psi} ,
      \bar{\Phi} ) \p \| g\|_2 \\
&\le& (1+\vare ) \p \bar{\Psi}^* (u^*) \p \delta (T|\Phi ,\Psi ).
\end{eqnarray*} \hz

To apply the results from section 1 we remark that for $T\in\comxy$

\[ \fty T \p = \p \sup_{\nu_{\Phi} (u:\lzn\nach X)=1} \ell (Tu)
   \plll \mbox{and} \plll
   \fco T \p = \p \sup_{\ell (u:\lzn\nach X)=1} \pif (Tu) \p . \]

Using $\pzn u) = \pi_{U_n} (u)$ and $t_2^n (T) = \delta (T|G_n,U_n) =
\delta (T^*|U_n,G_n)$ ($U_n$ and $G_n$ are defined in the introduction)
it is clear that (cf. \cite{TOM} (Theorem 25.5))

\[ t_2^n (T) \p = \p \sup_{ w\in\com Y^*),\ell^* (w^*) = 1} \pzn T^*w)
   \plll \mbox{and} \plll
   c_2^n (T) \p = \p \sup_{\ell (u:\lzn\nach X)=1} \pzn Tu) \p . \]

Furthermore, via $\fty T = \delta (T^*| \bar{\Phi},G_n)$ we obtain

\[ \ftyh T \p = \p \sup_{w\in\com Y^*),\ell^* (w^*) = 1}\pi_{\bar{\Phi}}
		(T^*w) \p . \]

Sometimes we will use

\[ \fty\len \kl \ftyh \len \kl K \sqrt{\log (n+1)} \p \fco \lin \]

which is an easy consequence of $\ftyh \len = \delta (\lin |\Phi ,G_n )$
and of $\ell (u) \le K \sqrt{\log (n+1)} \ell^*(u^*)$ in the case
$u\in\com \lin )$ (\cite{PS3}, \cite{TOM} (Theorem 12.7)). \hz

Let us start with the following standard lemma (cf. \cite{PIE}(6.2.7)).\hz

\begin{lemma} \label{phi_p2} Let $u\in\com X)$ and $\fis$ be an orthonormal
system. Then

\[ \pz u^*) \le \bar{\Phi}^* (u^*) \le \Phi (u) \le \pz u). \]

\end{lemma}

$Proof:$ For some normalized Borel measure $\mu$ on $\bzn$ (see
\cite{PIEOP}(17.3.2)) we get
\begin{eqnarray*}
\Phi (u)
& = & \kla \int_{\Om} \noo \summ_{k=1}^n \phi_k ue_k \rrm^2
      dp \mer^{\frac{1}{2}} \le \pz u) \kla \int_{\Om} \int_{\bzn}
      \left |\langle \summ_1^n \phi_k e_k ,a \rangle \right |^2 d\mu
      (a) dp \mer^{\frac{1}{2}} \\
&\le& \pz u) \kla \int \summ_1^n |\langle e_k ,a \rangle |^2 d\mu (a)
      \mer^{1/2}  \le \pz u).
\end{eqnarray*}
Using trace duality and $\pi_2^* (u^*) = \pz u^*)$ from
\cite{PIEOP}(19.2.14) we obtain $\pz u^*) \le \Phi^* (u^*)$ and in the same
way $\pz u^*) \le \bar{\Phi}^* (u^*)$ Finally, assuming $w\in\com X^*)$ the
inequality $\bar{\Phi}^* (u^*) \le \Phi (u)$ follows from
\for |{\rm tr}(u^*w)|\p =\p \left |\summ_{k=1}^n\langle ue_k ,we_k\rangle
     \right |
     \p =\p \left | \int \left\langle \summ_k \phi_k ue_k , \summ_l
     \bar{\phi}_l we_l \right\rangle dp \right |
     \p \le \p \Phi (u) \bar{\Phi} (w)\p .  \plll \Box \mel

Lemma \ref{phi_p2} together with $\pi_2(u) \kl \sqrt{2} \pzn u)$ for
$u\in\com X)$ (\cite{TOM}) imply the easy part of Theorem \ref{phi-type}.\hz

\begin{cor} \label{easy_estimate} For all $T\in {\cal L}(X,Y)$ one has
\[ \fty T \le \ftyh T \le \sqrt{2} t_2^n(T) \pll \mbox{and} \pll
   \fco T \le \sqrt{2} c_2^n(T). \]
\end{cor} \hz

We come to the non-trivial part.\hz

{\bf Proof of Theorem \ref{phi-type} in the introduction:}
(1) Since $\ftyh \len \ge \delta \sqrt{n}$ implies
$\fco \lin \ge \frac{\delta}{K} \sqrt{\frac{n}{\log (n+1)}}$
Theorem \ref{alpha_pz} (1) gives
\[ (\delta /K)^3 \pz w) \le c \pif (w) \hspace{.7cm}\mbox{for all}
   \hspace{.7cm} w\in \com X)\]
such that the conclusion with respect to the cotype will be clear.
In the type situation we have to observe that $\fco \lin = c_{\bar{\Phi}}
(\lin )$ and hence $(\delta /K)^3 \pz w) \le c \pi_{\bar{\Phi}} (w)$.
(2) One obtains $\delta^3 \nu_{\Phi} (w) \kl c \pz w^*)$
from Theorem \ref{alpha_pz} (2). Then we use for $T\in\comxy$ the equality
$t_2^n(T) \p=\p \sup\{ \ell(Tu) | \pi_2^n (u^*), \p u\in \com X)\}$,
see \cite{TOM} (Theorem 25.5,24.2) to conclude.
\hfill $\Box$ \hz

\begin{rem} \label{comp-cond}
Since $\fty\len \kl \ftyh \len \kl K \sqrt{\log (n+1)} \p \fco \lin$
the assumptions in (1) of
Theorem \ref{phi-type} are weaker than the assumption made in (2) of
Theorem \ref{phi-type}. Theorem \ref{sumfew} shows that the
assumptions of (1) are "strictly" weaker than the assumption of (2).
\end{rem} \hz

{\bf Proof of Theorem \ref{geometry} in the introduction:}
Let $w\in\com\lin )$ and
$h_j := \summ_{k=1}^n \langle we_k,e_j \rangle \phi_k$. Then
$\Phi (u) =\kla \sup_{j=1,...,n} |h_j (\om )|^2 dp(\om )\mer^{\frac{1}{2}}$
and $\|h_j\|_2 = \| w^*e_j\|$ such that
\begin{eqnarray*}
\pif (\lin )
& = & \sup \left\{ \Phi (w) \mitt w\in\com\lin ) , \| w\| \le 1 \right\} \\
& = & \sup \left\{ \kla \int_{\Om} \sup_{j=1,...,n} |h_j (\om )|^2 dp(\om )
      \mer^{\frac{1}{2}} \mitt (h_j)_1^n \subseteq {\rm span} \{ \phi_k \},
      \| h_j\|_2 \le 1 \right\}
\end{eqnarray*}
and
\begin{eqnarray*}
\pif (\izin )
& = & \sup \left\{ \Phi (\izin w) \mitt w\in\com\lzn ) , \| w\| \le 1 \right\}
\\
& = & \sup \left\{ \kla \int_{\Om} \sup_{j=1,...,n} |\psi_j (\om )|^2
      dp(\om ) \mer^{\frac{1}{2}} \mitt (\psi_j)_1^n \subseteq {\rm span}
      \{\phi_k\} \pl \mbox{orthonormal} \right\}.
\end{eqnarray*}
For the latter equality we use the fact that it is sufficient to take the
supremum over all orthogonal matrices $w\in\com\lzn )$. The implications
$(G_n) \rightarrow (G^{'}_n) \rightarrow (G^{''}_n)$ follow immediately from
Lemma \ref{l-norm_to_op-norm} and Corollary \ref{cor_delta^3} ($\pi_2 (\izin )
=\sqrt{n}$).
For $(G^{''}_n) \rightarrow (G_n)$ we use Theorem \ref{main_abs} to deduce
$\delta^{''} \pz u) \kl \sqrt{12} \p \pif (u)$ such that
$\delta^{''}  c_2^n (\lin ) \kl \sqrt{12} \p \fco \lin$. \hfill
$\Box$ \hz

Finally we prove the infinite versions of Theorem \ref{phi-type}. Before
doing this we need

\pagebreak

\begin{lemma} \label{appr} Let $(f_l)^n_1 \subset L_2(\Om,p)$ be a normalized
sequence and $H\subseteq L_2(\Om,p)$ be an $n$-dimensional subspace such that
\[\summ_1^n \bet \kla f_l, h \mer \rag^2  \pl \le \pl \noo h \rrm^2\]
for all $h \in H$. If $\sup\left\{ |(f_l,h)| \mitt h\in B_H\right\}\ge\theta$
for all $l=1,...,n$ then there exists an orthonormal basis  $(\psi_l)_1^n$ of
$H$ and
a subset $I\subseteq\{1,...,n\}$ with $|I|\ge \frac{\theta^6}{c} \p n$
such that
\[ |(f_l,\psi_l)| \p \ge \p \frac{\theta^3}{c} \hspace{1.1cm} \mbox{for all}
   \hspace{1.1cm} k\in I \p ,\]
where $c\ge 1$ is an absolute constant.
\end{lemma} \hz

$Proof:$ Our assumption ensures $(f_l, h_l) \ge \theta$ for some
$h_1,...,h_n \in B_H$. We fix an isometry $T:\ell_2^n \nach H$ and define a
norm $\beta$ on
${\cal L}(\ell_1^n,\ell_2^n)$ by
\[ \beta(u ) \pl :=\pl \sup_{\noo w:\ell_2^n \nach \ell_2^n \rrm \le 1}
   \kla \summ_1^n \bet \kla f_l, Twu(e_l) \mer \rag^2 \mer^{\frac{1}{2}}
   \pl .\]
Let us note that the first inequality of our assumption gives
$\beta(a\otimes x)\le \noo a \rrm_{\lin}\noo x \rrm_{\lzn}$, which implies
$\beta \le \nu$ on the component ${\cal L}(\ell_1^n,\ell_2^n)$. Furthermore,
the operator
$u\p:=\p \summ_1^n  e_l \otimes T^{-1}(h_l)\in
{\cal L}(\len,\lzn)$ is of norm at most one and satisfies
$\beta(u)\ge \p\theta \p \sqrt{n}$. In this situation we can apply Lemma
\ref{lemma_delta^3} to deduce for some $c\ge 1$
\[\beta(\iezn)\pl\ge\pl \frac{\theta^3}{c}\p \sqrt{n} \pl .\]
By convexity we find an orthogonal matrix $O$ such that
$\summ_1^n \bet\kla f_l,TO(e_l) \mer \rag^2 \pl \ge \pl \frac{\theta^6}{c^2}
 \pl n$ . Clearly,
\[ \psi_l\p:=\p TO(e_l) \p \frac{\bet\kla TO(e_l) , f_l \mer \rag}
   {\kla TO(e_l) , f_l \mer} \p , \]
where we assume $\frac{0}{0} = 1$ , defines an orthonormal basis in $H$. From
$\bet\kla f_l,TO(e_l) \mer \rag\le 1$ we derive that the set
\[ I\pl := \p \left \{ l \mitt \pl \kla \psi_l , f_l \mer
   \p\ge\p \frac{\theta^3}{c\p \sqrt{2}} \right \}\]
is of cardinality at least $\frac{\theta^6}{2\p c^2} \p n $. \hfill $\Box$ \hz

\pagebreak


\begin{theorem} \label{cotype_infinite} Let $\Phi = (\phi_k )_{k\in\nz}
\subset L_2(\Om ,p)$ be an orthonormal system. Then the following assertions
are equivalent.
\begin{enumerate}
\item [(1)] There exists a constant $c>0$ such that $c_2 (T) \le c \p \fco T$
      for all operators $T\in\comxy$.
\item [(2)] There exists $\delta >0$ such that $\fco \lin \ge \delta
      \sqrt{\frac{n} {\log (n+1)}}$ for all $n=1,2,...$
\item [(3)] There exists $\eta >0$ such that for all $n=1,2,...$ there is an
      orthonormal system $\Psi = (\psi_l )_{l=1}^n \subset
      {\rm span} \{\phi_k \}$ and disjoint measurable subsets $A_1,...,A_n$
      such that
      \[\int_{A_l} |\psi_l|^2 dp \ge \eta \pla \mbox{for} \pll l=1,...,n \p
      .\]
\item [(4)] There exists $\theta >0$ such that for all $n=1,2,...$ there is an
      $n$-dimensional subspace $H \subseteq {\rm span} \{\phi_k \}$ and
      an orthonormal system $f=(f_l)_1^n \subset L_2(\Om,p)$ such that the
      $f_l$ have disjoint support and
      \[ \sup_{h\in B_H} |(f_l,h)| \p \ge \p \theta \pla \mbox{for} \pll
	 l=1,...,n \p .\]
\end{enumerate}
\end{theorem}\hz

$Proof:$ $(1) \rightarrow (2)$ is trivial. $(2)\rightarrow (3)$ For fix
$n$ there are $x_1,...,x_N \in \lin$ such that

\[ \noo\summ_1^N \phi_k x_k \rrm_2 \ge \frac{\delta}{2}
   \sqrt{\frac{n}{\log (n+1)}} \pll \mbox{and} \pll
   \noo \summ_1^N g_k x_k\rrm_2 \le 1 \p .\]

Assuming ${\rm span}\{ x_1,...,x_N\}=\lin$ it is easy to see that there are
$y_1,...,y_n\in\lin$ and a matrix $(p_{ij})_{i=1,j=1}^{n,N}$ such that
$\summ_{j=1}^N p_{ij} \bar{p_{kj}} = \delta_{ik}$ and
$x_j = \summ_i p_{ij} y_i$ for $j=1,...,N$. Consequently,

\[ \noo\summ_1^n \psi_i y_i \rrm_2 \p = \p \noo\summ_1^N \phi_j x_j \rrm_2
   \p \ge \p \frac{\delta}{2} \sqrt{\frac{n}{\log (n+1)}}
   \pll \mbox{and} \pll
   \noo\summ_1^n g_i^{'} y_i\rrm_2\p =\p \noo \summ_1^N g_j x_j\rrm_2\p\le
   \p 1 \]

if $\Psi := (\psi_i)_1^n = \kla \summ_{j=1}^N p_{ij} \phi_j\mer_{i=1}^n$.
Applying Theorem \ref{alpha_pz} yields
$\kla \frac{\delta}{2} \mer^3 \pz w) \le c \pi_{\Psi} (w)$ for all
$w\in\com X)$. Especially,
$\pi_{\Psi} (\izin ) \ge \frac{1}{c} \kla\frac{\delta}{2}\mer^3 \sqrt{n}$
such that there is an orthonormal system $(h_k)_1^n \subset {\rm span}
\{\psi_i \}$ with

\[ \int \sup_k |h_k|^2 dp \p \ge \p \kla\frac{1}{c}\mer^2
   \p \kla \frac{\delta}{2}\mer^6 n \p . \]

Applying \cite{TOM} (Lemma 31.3) we find an index set $J\subseteq\{ 1,...,n\}$
with $|J|\ge \frac{(\delta /2)^6}{2c^2} n$ and disjoint measurable sets $A_k$
such that $\int_{A_k} |h_k|^2 dp \ge \frac{(\delta /2)^6}{2c^2}$ for $k\in
J$.
$(3)\rightarrow (1)$ It is clear that we have

\[ \int \sup_k |\psi_k |^2 dp \p \ge \p \eta n \]

such that $\pi_{\Psi}(\izin ) \ge \sqrt{\eta} \sqrt{n}$. Applying  Theorem
\ref{main_abs} we obtain $\pz u) \eta^{\frac{1}{2}} \le c \pi_{\Psi} (u)$
for all $u\in\com Y)$. Assuming $\psi_i = \summ_{j=1}^N p_{ij} \phi_j$ such
that $\summ_{j=1}^N p_{ij} \bar{p}_{kj} = \delta_{ik}$ we obtain for
$P:=\kla p_{ij} \mer\in{\cal L}(\ell_2^N ,\lzn )$, $u\in\com X)$,
$T\in\comxy$ and some $A\in\com\lzn )$ with $\| A\| =1$
\[ \pz Tu) \kl c \p \eta^{-\frac{1}{2}} \p \pi_{\Psi} (Tu) \p
   =  \p c \p \eta^{-\frac{1}{2}} \p \Psi (TuA) \p
   =  \p c \p \eta^{-\frac{1}{2}} \p \Phi (TuAP) \]
such that
\for \pz Tu) \kl c \p \eta^{-\frac{1}{2}} \p \fco T \p \ell (uAP) \kl c \p
     \eta^{-\frac{1}{2}} \p \fco T \p \ell (u) \p.
\mel
$(3)\rightarrow (4)$ We take $H:= {\rm span} \left\{ \psi_l \mitt l=1,...,n
\right\}$ and $f_l := \frac{\psi_l \chi_{A_l}}{\| \psi_l \chi_{A_l} \|}$
such that $|(f_l ,\psi_l)| \ge \sqrt{\eta}$.
$(4)\rightarrow (3)$ We apply Lemma \ref{appr} and get an orthonormal basis
$\Psi = (\psi_l)_1^n$ of $H$ and a proportional subset $I\subseteq \left\{
1,...,n\right\}$ with $|I| \ge \frac{\theta^6}{c} n$ such that for $l\in I$
and $A_l := {\rm supp} (f_l)$ one obtains
\[ \frac{\theta^3}{c} \kl (f_l,\psi_l) \kl \| f_l\| \p \| \psi_l \chi_{A_l} \|
   \p .\]
Hence $\int_{A_l} |\psi_l|^2 dp \p \ge \p \frac{\theta^6}{c^2}$ for $l\in I$.
\hfill $\Box$ \hz

\begin{theorem} \label{type_infinite} Let $\Phi = (\phi_k )_{k\in\nz}
\subset L_2(\Om ,p)$ be an orthonormal system. Then the following assertions
are equivalent.
\begin{enumerate}
\item [(1)] There exists a constant $c>0$ such that
      $t_2 (T) \le c \p t_{\Phi}(T)$ for all operators $T\in\comxy$.
\item [(2)] There exists $\delta >0$ such that $t_{\Phi}(\len )\ge \delta
      \sqrt{n}$ for all $n=1,2,...$
\end{enumerate}
\end{theorem}\hz

$Proof:$ Clearly it remains to show $(2)\rightarrow (1)$. Using the argument
as in $(2)\rightarrow (3)$ of the above theorem we find an orthonormal system
$\Psi_n = (\psi_i)_1^n \subset span\{\phi_j\}$ and $y_1,...,y_n\in\len$ such
that
\[ \noo\summ_1^n \psi_i y_i \rrm_2 \kl 1
   \hspace{.8cm} \mbox{and} \hspace{.8cm}
   \noo\summ_1^n g_i^{'} y_i\rrm_2\p \ge \p \frac{\delta}{2} \sqrt{n} . \]
Theorem \ref{alpha_pz} yields
\[ \kla \frac{\delta}{2} \mer^3 \nu_{\Psi_n} (w) \le c \pi_2(w^*)
   \hspace{.7cm}\mbox{for all} \hspace{.7cm} w\in\com X).\]
Consequently, for all $T\in\comxy$ and all $n=1,2,...$
\for \frac{\delta^3}{8c} t_2^n (T) \kl t_{\Psi_n} (T) \kl t_{\Phi} (T) .
     \\[-1.5cm] \mel \hfill $\Box$ \hz

\begin{rem} \label{nontrivial_ann} In Proposition \ref{nontrivial_sys} we will
see that there exists an orthonormal system
$\Phi = (\phi_k)_1^{\infty} \subseteq (e^{ikt})_{k\in\nz} \subseteq L_2(\Pi )$
such that the $\Phi$- cotype and type does not coincide with the usual
cotype 2 and type 2, but is nontrivial, i.e. there are Banach spaces without
$\Phi$- cotype and type.
\end{rem}\hz


\section{Orthonormal systems on discrete measure spaces}\label{few}

\setcounter{lemma}{0}

In this section we compare (sometimes for simplicity in the real situation)
ordinary type and cotype constants with the $\Phi$-type and cotype constants
in the case that the orthonormal system $\Phi$ lives on a discrete measure
space $\Om := \{ \om_1 ,...,\om_N\}$. We start with the positive part by
showing that the $\pi_{\Phi}$-norm and the $\pi_2^n$-norm are close to each
other whenever $n\sim N$ (which clearly implies the same for the corresponding
cotype constants). In order to apply Theorem \ref{main_abs} we need the
following lemma which contains an argument discovered in a discussion
with B. Kashin. \hz

\begin{lemma} \label{induct} Let $1\le n \le N$ and let $H\subset \ell_2^N$ be
an $n$-dimensional subspace. Then there exists an orthonormal basis
$(h_k)_1^n$ and pair wise different coordinates $j_k \in \{1,..., N\}$ such
that
\[ \bet \kla h_k ,e_{j_k}\mer \rag^2 \pl \ge \pl \frac{n}{3N} \pll
   \mbox{for all } \pll k\p\le \p \frac{n}{3} \pl .\]
\end{lemma}

$Proof:$ Let $P_H$ be the orthogonal projection onto $H$. It is well-known
that
\[ n \pl =\pl \pi_2(P_H)^2 \pl =\pl \summ_1^N \noo P_H(e_j) \rrm^2 \pl .\]
Hence there exists $j_1\in \{ 1,...,N\}$ such that
\[ \noo P_H(e_{j_1}) \rrm^2 \pl \ge \pl \frac{n}{N} \pl .\]
We consider the normalized element
$h_1\:=\p \frac{P_H(e_{j_1})}{\noo P_H(e_{j_1})\rrm}$
which satisfies ($P_H$ is a projection)
\[ \kla h_1 ,e_{j_1} \mer \pl = \pl
   \frac{\kla P_H(e_{j_1}),e_{j_1}\mer}{\noo P_H(e_{j_1})\rrm }
   \pl =\pl  \noo P_H(e_{j_1}) \rrm \pl \ge \pl
   \sqrt{\frac{ {\rm dim}(H)}{N}} \pl . \]
Now we can proceed by induction setting $H^1 :=H$ and
$H^{k+1}:=H^{k}\cap {\rm span}\{ e_{j_k},h_k\}^{\bot}$. Here $j_k$ and
$h_k\in H_j$ are chosen by the construction above and satisfy
$\noo h_k \rrm \p=\p 1$ and
\[ |(h_k ,e_{j_k} )|^2 \pl \ge \pl \frac{{\rm dim}(H_k)}{N}
   \pl \ge \pl \frac{n-2k+2}{N} \pl .\]
If we continue as far as $k\le\frac{n}{3}$ we get an orthonormal sequence
$(h_k)_{k\le\frac{n}{3}}$ in $H$ with the desired properties. Note that by
construction the elements $e_{j_k}$ are disjoint. Finally we complete the
sequence $(h_k)$ to an orthonormal basis of $H$. \hfill $\Box$ \hz

Now we compare the $\pi_{\Phi}$-norm with the $\pi_2^n$-norm. \hz

\begin{prop} \label{p2_few} Let $\fis\subset \ell_2^N$ be an orthonormal
system. Then
\[ \pi_2^n(u) \kl 12 \sqrt{\frac{N}{n}} \p\p\p \pif (u)
   \hspace{.7cm} \mbox{for all} \hspace{.7cm} u\in\com X) \p .\]
\end{prop} \hz

$Proof:$ First we show $\sqrt{n} \pl \le \pl 2\sqrt{3} \pl\sqrt{\frac{N}{n}}
\pl \pi_{\phi}(\iota_{2,\infty}^n)$ .
Since $\pi_{\phi}(\iota_{2,\infty}^n) \ge 1$ we can assume $n\ge 12$.
Setting $H \p:=\p {\rm span}\{ \phi_k \}$ we choose an orthonormal system
$(h_k)_1^n$ and pair wise different $k_j$ according to Lemma \ref{induct} .
Now the lower estimate of $\pi_{\phi}(\iota_{2,\infty}^n)$ follows from

\for
\pi_{\phi}(\iota_{2,\infty}^n)^2
& = & \sup \left\{ \summ_{j=1}^N \sup_{k=1,...,n} |\psi_k (j)|^2
      \mitt (\psi_k)_1^n \subseteq {\rm span} \{\phi_k\}_1^n \pl
      \mbox{orthonormal} \right\} \\
&\ge& \summ_{j=1}^N \sup_{k=1,..,n} \bet (h_k,e_j) \rag^2
 \ge  \summ_{k \le \frac{n}{3}} \bet (h_k,e_{j_k}) \rag^2
 \ge  \frac{n}{4} \pl \frac{n}{3N} \p .\\
\mel

Finally, Theorem \ref{main_abs} yields the desired assertion.
\hfill $\Box$ \hz

The above proposition gives
\[ c_2^n (T) \kl 12 \sqrt{\frac{N}{n}} \pif (T) \pll \mbox{whenever} \pll
   \fis\subset\ell_2^N \]
which was claimed in Theorem \ref{sumfew}(1). To prove the remaining
parts of Theorem \ref{main_abs} we have to construct orthonormal systems with
small type or cotype constants.
The notion of a $\La_p$-system originally introduced for Fourier series
turns out to be a useful tool. For $1 \! \le\! p\! <\! \infty$ an
orthonormal system $\Phi = (\phi_i )_{i\in I}$ ($I$ is a countable index set)
on a probability space $(\Om,p)$ is said to be
a $\La_p$-system if there exists a constant $c\!\ge\!0$ such that for all
finitely supported sequences $(\al_i)_I \subset \kz^I$ one has

\[ \kla \int\limits_{\Om} \bet \summ_i \al_i \phi_i \rag^p dp \mer^{1/p}
 \le \pl c \pl \int\limits_{\Om} \bet \summ_i \al_i \phi_i \rag dp \]

By $\La_p(\Phi)$ we denote the best constant in the inequality above.
In order to construct orthonormal systems with small cotype constants we also
need the notion of a $K_q$-system. For $2 \! \le\! q\! <\! \infty$ an
orthonormal system $\Phi = (\phi_i )_{i\in I}$ on a probability space
$(\Om,p)$ is said to be a $K_q$-system if there exists a constant $c\!\ge\!0$
such that for all finite sequences $(\al_i)_I \subset \kz^I$ one has

\[ \kla \int\limits_{\Om} \bet \summ_i \al_i \phi_i \rag^q d\mu \mer^{1/q}
 \le \pl c \pl \kla \summ_i \bet \al_i \rag^2 \mer^{1/2} . \]

By $K_q(\Phi)$ we again denote the best constant in the inequality above.
In fact for $q\!>\!2$ every $K_q$-system is a $\La_q$-system and vice versa.
Note that any finite orthonormal system, that is the index set $I$ is assumed
to be finite, $\Phi\subset L_p$ is a $\Lambda_p$-
system and any finite orthonormal system $\Phi\subset L_q$ is a $K_q$- system
but the constants could be different and are important in the sequel.
In the following it will be convenient to use
$L_p^N := [\kz^N , \| \pl \|_{L_p^N}]$ with

\[ \| (\xi_k )\|_{L_p^N} :=\kla \frac{1}{N} \summ_1^N |\xi_k |^p
    \mer^{\frac{1}{p}} \]

instead of $\ell_p^n$. We will start with the construction of orthonormal
systems
$\Phi$ such that $\fco \lin$ is small.


\begin{lemma} \label{Kaq} Let $2\!<\!q\!<\!\infty$ and let
$\Phi = (\phi_k)_1^N$ be an orthonormal system. Then
\[ \fco \lin \kl c_q \p \frac{n^{\frac{1}{q}}}{\sqrt{\log (n+1)}} \p
    K_q(\Phi) \p,\]
where $c_q>0$ is an absolute constant depending on $q$ only.
\end{lemma} \hz

$Proof:$ (1) First we show for $u\in\com\lin )$
\[ \pi_q (u) \kl c \p \frac{n^{\frac{1}{q}}}{\sqrt{\log (n+1)}} \p \ell (u)
   \p .\]

Applying \cite{LE-TA}(Theorem 12.10) in the situation

\[ X_t \p := \p \left\langle \summ_{i=1}^n g_i ue_i ,e_t \right\rangle
   \pll \mbox{for} \pll t\in T:=\{ 1,...,n\} \p ,\]

where $(e_t )_1^n$ is the unit vector basis of $\len$, yields a sequence
$(Y_k)_{k\ge 1}$ of gaussian variables with $\| Y_k\|_2 \le
\frac{c \p\ell (u)}{\sqrt{\log (k+1)}}$ such that

\[ X_t \p \stackrel{L_2}{=} \p \summ_{k\ge 1} \p \al_k(t) \p Y_k \]

for $\al_k(t)\ge 0$ with $\summ_k \al_k(t) \le 1$ (in the complex case we
consider the real and complex part separately and obtain  complex
$\al_k(t)$ with $\summ_k |\al_k(t)|\le 1$). Setting
$u_k := \sqrt{\log (k+1)} (\langle Y_k,g_i\rangle )_{i=1}^n \in \lzn$,
$v_k := (\al_k(t))_{t=1}^n \in\lin$,
$   A:= \summ_{k\ge 1} u_k \otimes e_k \in\com \ell_{\infty} )$,
$   D:= \summ_{k\ge 1} \frac{1}{\sqrt{\log (k+1)}} e_k \otimes e_k \in
	{\cal L} (\ell_{\infty} ,\ell_{\infty} )$,
$   B:= \summ_{k\ge 1} e_k \otimes v_k \in {\cal L}(\ell_{\infty} ,\lin )$
we deduce $\| u_k\| \le \sqrt{\log (k+1)} \| Y_k\|_2 \le c \p\ell (u)$,
$\| A\| \le c \ell (u)$, $\| B\| \le 1$, and a factorization

\[ u\p : \p \lzn \nachop A \ell_{\infty} \nachop D
   \ell_{\infty} \nachop B \lin \p .\]

Considering $D=D_1 + D_2$, whereas $D_1$ is the diagonal operator associated
to the sequence $\kla \frac{1}{\sqrt{\log 2}},...,
\frac{1}{\sqrt{\log (n+1)}},0,...\mer$ we obtain from \cite{DE-JU}(Theorem 5)
and $\pi_q (D_a) \kl \| a\|_q $, if the diagonal operator $D_a \in
{\cal L}(\ell_{\infty},\ell_{\infty})$ is generated by the sequence $a$,

\begin{eqnarray*}
\pi_q (u)
&\le& \| B\| \p \pi_q (D_1) \p \| A\| \p + \p \pi_q^n (B) \p
      \frac{1}{\sqrt{\log (n+1)}} \p \| A\| \\
&\le& c_q \p \| B\| \p \frac{n^{\frac{1}{q}}}{\sqrt{\log (n+1)}} \p \| A\| \p
      + \p n^{\frac{1}{q}} \p \| B\| \p
      \frac{1}{\sqrt{\log (n+1)}} \p \| A\| \\
&\le& c^{'}_q \frac{n^{\frac{1}{q}}}{\sqrt{\log (n+1)}} \p \ell (u)\p .
\end{eqnarray*}
(2) Assuming $v\in{\cal L}(\ell_2^N,\lin )$ it is easy to see that
there is a factorization $v=uP$ where $P\in{\cal L}(\ell_2^N ,\lzn )$ and
$u\in\com \lin )$ such that $\| P\| = 1$ and $\ell (u) =\ell (v)$.
Using the continuous version of the q-summing norm (which follows by
an easy approximation argument, see \cite{PIS}(Proposition 1.2)) we can
deduce

\begin{eqnarray*}
\Phi (v)
&\le& \kla \int\limits_{\Om} \noo \summ_1^N \phi_k v(e_k) \rrm^q dp \mer^{1/q}
      \kl \pi_q(v) \sup_{\noo (\al_k) \rrm_2 \le 1} \kla \int\limits_{\Om}
      \bet \summ_1^N \al_k \phi_k \rag^q dp \mer^{1/q} \\
&\le& \pi_q(v) \pl K_q(\Phi ) \kl \pi_q(u) \pl K_q(\Phi ) \kl  c^{'}_q  \p
      \frac{n^{\frac{1}{q}}}{\sqrt{\log (n+1)}} \p K_q (\Phi ) \ell (u) \\
&\le& c^{'}_q \p\frac{n^{\frac{1}{q}}}{\sqrt{\log (n+1)}} \p K_q (\Phi )
      \ell (v) \pl .\\[-1.5cm]
\end{eqnarray*}
\hfill $\Box$ \hz

The following proposition provides small orthonormal systems with small
cotype constants. \hz


\begin{prop} \label{smallcotype} Let $2\!<\!q\!<\!\infty$ and $1\!\le\!n\!\le
\!N$. Then there exists an real orthonormal system $\fis\subset
L_2^N$ with
\begin{enumerate}
\item [(1)] $\frac{\sqrt{n}}{N^{1/q}} \kl K_q(\Phi) \le \pl c \pl
      \max \{\sqrt{q}, \frac{\sqrt{n}}{N^{1/q}} \}, $
\item [(2)] $\fco \lin \sqrt{\frac{\log(n+1)}{n}} \kl c_q \pl \max
      \left\{\sqrt{q}\, n^{1/q-1/2}, \kla \frac{n}{N} \mer^{1/q} \right\}, $
\end{enumerate}

where $c>0$ is an absolute constant and $c_q>0$ depends on $q$ only.
In particular, for $0<\delta<1$ and $N= [n^{1/\delta}]$ one has
$\fco \lin \kl c_{\delta} \p n^{\frac{\delta}{2}}$\p .
\end{prop} \hz

$Proof:$ (1) We will use a random argument. By the comparison principle
for random orthonormal matrices and gaussian variables of Marcus and Pisier,
see \cite{B-G}, and Chevet's inequality, see\cite{CHV},  we deduce
\begin{eqnarray*}
 \sqrt{N} \, \ew \noo \summ_{k=1}^n \summ_{j=1}^N o_{k,j} e_k \otimes
       e_j\::\:\lzn \to L_q^N \rrm
 &\le& \frac{c_0}{N^{\frac{1}{q}}}\, \ew \noo \summ_{k=1}^n \summ_{j=1}^N
       g_{k,j} e_k \otimes e_j\::\:\lzn \to \ell_q^N \rrm \\
 &\le& \frac{c_0\pl \sqrt{2}}{N^{\frac{1}{q}}} \, \kla \ew \noo \summ_{j=1}^N
       g_j e_j \rrm_{\ell_q^N} \:+\: \ew \noo \summ_{k=1}^n g_k e_k
       \rrm_{\ell_2^n} \mer\\
 &\le& c_0\pl c_1\pl \sqrt{2}\pl \kla \sqrt{q}\:+\:\sqrt{n}\p N^{-\frac{1}{q}}
 \mer \,.\\
\end{eqnarray*}
Here the expectation is taken with respect to the Haar-measure on the group
${\cal O}(N)$ of orthonormal matrices and with respect to the standard
gaussian density in $\rz^{nN}$. For a random matrix $o$ satisfying the above
inequality we define the orthonormal system $\Phi = (\phi_k)_1^n \subset
L_2^N$ by
\[ \phi_k\,:=\,\sqrt{N} \, \summ_{j=1}^N o_{k,j} \,e_j \p . \]

Therefore we have proved
\[ \ew \p K_q(\Phi) \,\le \, 3 \pl c_0 \pl c_1\pl \max\{\sqrt{q},\sqrt{n}\,
   N^{-1/q} \} \pl .\]

and can choose our system randomly (with an obvious change of constants the
same estimate for the $K_q$-constant is valid if we compute this constant
with complex coefficients). To obtain the lower estimate of the
$K_q$-constant we first claim that $\sqrt{n} \kl \pif (\ell_{\infty}^N)$.
In order to prove this claim we find $(a_j)_1^N \subset \bzn$ such that
$\left\langle \ssum_{k=1}^n\phi_k(j) e_k, a_j \right\rangle \p=\p \noo
 \ssum_{k=1}^n \phi_k(j) e_k \rrm_{\lzn} $
and define the operator
$R\p:=\p \ssum_1^N a_j \otimes e_j \in \com \ell_{\infty}^N)$ with
$\| R\| \le 1$. Then it follows that

\[ \sqrt{n} \pl =\pl \kla \ssum_1^n \noo \phi_k \rrm_{L_2^N}^2 \mer^
   {\frac{1}{2}}
   \pl = \pl
   \kla \frac{1}{N} \ssum_{j=1}^N \noo \ssum_{k=1}^n \phi_k(j)
   e_k\rrm_{\ell_2^n}^2 \mer^
   {\frac{1}{2}}
   \kl \Phi(R) \kl \pif (\ell_{\infty}^N) \p . \]

Using the argument given in the end of the proof of Lemma \ref{Kaq} we
continue with

\[ \pif (\ell_{\infty}^N) = \sup \left\{ \Phi (u) \p\mitt\p \noo u:\lzn
   \rightarrow \ell_{\infty}^N\rrm \le 1 \right\} \le
   K_q(\Phi ) \sup \left\{ \pi_q (u) \p\mitt\p \noo u:\lzn
   \rightarrow \ell_{\infty}^N \rrm \le 1 \right\}  . \]

Finally, \cite{DE-JU}(Theorem 5) implies $\pi_q (u) \kl N^{1/q} \| u\|$.
(2) is a consequence of (1) and Lemma \ref{Kaq}. The last assertion follows
by $q=\frac{2}{\delta}$. \hfill $\Box$\hz

We continue by constructing orthonormal systems with
small type constants in the proportional case. \hz

\begin{lemma} \label{lam2} Let $\fis$ be an orthonormal system. Then
$\fty {\ell_1} \kl c \p \Lambda_2 (\Phi )$ for some absolute
constant $c>0$.
\end{lemma}\hz

$Proof:$ Let $x_1,...,x_n \subset \ell_1$. Using the Kahane
inequality for gaussian averages due to Hoffmann-J\H{o}rgensen, see
\cite{LE-TA}, we deduce
\for
\noo \summ_1^n g_k x_k \rrm_{L_2(\ell_1)} &\le& c \int \summ_{j \in \nz}
\bet \summ_{k=1}^n g_k \langle x_k,e_j \rangle \rag dp
\kl c \p \summ_{j \in \nz} \kla \summ_{k=1}^n \bet \langle x_k,e_j \rangle
\rag^2
      \mer^{\frac{1}{2}}\\
&\le& c \p \La_2(\Phi) \pl \summ_{j \in \nz} \int\limits_{\Om} \bet
      \summ_{k=1}^n \phi_k \langle x_k,e_j \rangle \rag dp \p
=\p   c \p \La_2(\Phi) \fme \pl.\\[-1.5cm]
\mel \hfill $\Box$ \hz

A more abstract version of this argument can be applied for Banach lattices
with finite cotype, see \cite{JU}. In contrast to the previous
results large orthonormal systems will now be constructed in $L_2^N$ .

\begin{prop} \label{smalltype} For all $1\!\le\!n\!<\!N$ there exists
an real orthonormal system $\fis\subset L_2^N$ such that
for some absolute constant $c>0$

\begin{enumerate}
\item [(1)] $ \La_2(\Phi) \le \pl c \pl \sqrt{ \frac{N}{N-n} \log
      (1+\frac{N}{N-n} ) } \pl,$
\item [(2)] $ \fty {\ell_1} \le \pl c \pl \sqrt{ \frac{N}{N-n} \log
      (1+\frac{N}{N-n} ) } \pl .$
\end{enumerate}

In particular, for $0\! < \! \vare \! < \! 1$ and $1\!\le\! n \! \le \!
(1-\vare)N$ there is an orthonormal system satisfying

\[ \fty {\ell_1} \le \pl c \pl \sqrt{\frac{1}{\vare} \p
   \log \kla 1+\frac{1}{\vare} \mer }\pl.\]
\end{prop} \hz

$Proof:$ By \cite{CP}(Theorem 2.2) there exists a subspace
$E \subset \ell_1^N $ with $n=N-m =  {\rm dim}(E) $ such that for all $x \in
E$

\[ \noo x \rrm_2 \le \pl c \pl \sqrt{ \frac{\log \kla 1 + \frac{N}{m} \mer
}{m} }
   \pl \noo x \rrm_1 \pl .\]

Now we consider E as a subspace of $L_1^N$, $L_2^N$, respectively. Then we
have for all $x \in E$

\[\noo x \rrm_{L_2^N} \pl =\pl \frac{1}{\sqrt{N}} \pl \noo x \rrm_2
  \kl \frac{c}{\sqrt{N}} \pl \sqrt{ \frac{\log\kla 1 + \frac{N}{m}\mer}{m} }
  \pl \noo x \rrm_1
  \pl =\pl c \pl \sqrt{ \frac{N}{m} \log \kla 1 + \frac{N}{m}\mer }
  \pl \noo x \rrm_{L_1^N} \pl . \]

If we choose an orthonormal system $\fis$ in E with respect to the
scalar product of $L_2^N$ we obtain for all sequences $(\al_k)_1^n \in \rz^n$
(and with the constant $2c$ instead of $c$ also for complex $(\al_k)$)
\[\kla \summ_1^n \bet \al_k \rag^2 \mer^{1/2} \pl =\pl
		 \noo \summ_1^n \al_k \phi_k \rrm_{L_2^N}
  \kl c \pl \sqrt{ \frac{N}{m} \log \kla 1 + \frac{N}{m} \mer } \pl
     \noo \sum_1^n\al_k \phi_k \rrm_{L_1^N} \pl . \]

Therefore assertion (1) is proved. Assertion (2) follows immediately from
Lemma \ref{lam2} . \hfill $\Box$ \hz

For $n = \delta N$ we can again choose a random orthonormal
system in $L_2^N$ satisfying the assertion of the above proposition, because
for random subspaces the corresponding norm estimate is valid, see \cite{Sz}
(cf.\cite{PS2}(Theorem 6.1)). We are now in position to complete the \hz

{\bf Proof of Theorem \ref{sumfew} in the introduction:}
(1) follows from Proposition \ref{p2_few} . (2) and (3) are consequences of
Proposition \ref{smallcotype} (2) and \ref{smalltype}. \hfill $\Box$ \hz

Analyzing Theorem \ref{phi-type} one can ask whether it is possible or not to
replace the space $\lin$ or $\len$ by an arbitrary space. We will show in
Proposition \ref{pos_results} that this is not possible in some sense.
Let us begin with the following construction which yields $n$-dimensional
quotients of $\ell_1^N$ with large $\Phi$-type constants. Given an arbitrary
orthonormal system $\fis \subset L_2(\Om,p)$ we consider the convex body

\[ \bfi \pl :=\pl \overline{\rm absconv}\left\{
   \frac{\sum\limits_1^n \phi_k(\om) e_k}
   {\kla \sum\limits_1^n\bet\phi_k(\om) \rag^2 \mer^{ \frac{1}{2} } }\pl
   \mitt\pl \summ\limits_1^n\bet\phi_k(\om) \rag^2 >0\right\} \pl \subset
   \pl \lzn \]

and the associated Banach space $E_{\Phi} \p :=\p [\kz^n,B_{\Phi}]$ via
the Minkowski functional of $B_{\Phi}$ (we will see that $E_{\Phi}$ is
correctly defined). From the definition it is clear that, for example, $\bfi$
is the convex combination of at most $N$ points if the measure space has $N$
atoms and no continuous part. The next lemma summarizes some properties of
the convex body $\bfi$.\hz
\begin{lemma} \label{phi_to_mu} Let $\fis$ be an orthonormal system. Then the
following holds true.
\begin{enumerate}
\item [(1)] There exists a normalized Borel measure $\mu$ supported by
      $S_{n-1} \cap B_{\Phi}$ satisfying \linebreak
      $\int_{\bzn} \langle x,e_i \rangle \overline{\langle x,e_j \rangle} d\mu
      (x)
      = \frac{1}{n} \delta_{ij}$ and, for all Banach spaces $X$ and all
      $u\in\com X)$,
      \[ \Phi (u) = \sqrt{n} \kla \int_{\bzn} \| ux\|^2 d\mu (x) \mer^{1/2}.\]
\item [(2)] If $\iota_{\Phi}\p:\p\lzn \nach \efi$ is the formal identity then
      $\| \iota_{\Phi}^{-1}\| \le 1$ and $\pi_2(\iota_{\Phi}^*)\le\sqrt{n}$.
\item [(3)] In the real situation one has
      \[ \frac{\ell (\iota_{\Phi})}{\sqrt{n}} \gl  \kla \frac{{\rm vol}(\bzn)}
	 {{\rm vol}(B_{\Phi})} \mer^{1/n} \quad \mbox{and} \quad
	 \frac{r(\iota_{\Phi})}{\sqrt{n}} \gl\frac{1}{c}\pl \kla
	 \frac{{\rm vol}(\bzn)}{{\rm vol}(B_{\Phi})}
	 \mer^{1/n} \pl ,\]
where $c>0$ is an absolute constant,
$r(u) := \kla \ew \noo \summ_1^n \vare_k ue_k \rrm^2 \mer^{\frac{1}{2}}$, and
$(\vare_k)_1^n$ is the sequence of independent random variables with
$ p(\vare_k=1)=p(\vare_k=-1)=\frac{1}{2}$.
\end{enumerate}
\end{lemma}\hz

$Proof:$ (1) Considering $F : [\Om ,{\cal F},p] \rightarrow
[\lzn ,{\cal B}(\lzn )]$ with $F(\om ):= (\phi_1 (\om ),...,\phi_n (\om ))$
we obtain an image measure $p_F$ on ${\cal B}(\lzn )$ such that
\[ \int_{\lzn} \frac{\| x\|^2_2}{n} dp_F (x) =
   \frac{1}{n} \int_{\Om} \summ_1^n |\phi_k (\om )|^2 dp(\om ) =1. \]
Hence $d\nu (x) := \frac{\| x\|^2_2}{n} dp_F (x)$ defines a normalized Borel
measure on $\lzn$ with $\nu (\{ 0\} )=0$. Finally, the image measure
$\mu$ of $\nu$ with respect to
\[ P:[\lzn \backslash \{ 0\} ,{\cal B} (\lzn\backslash\{ 0\} )] \rightarrow
   [\bzn ,{\cal B} (\bzn )] \hspace{1cm} \mbox{whereas} \hspace{1cm}
   P(x) := \frac{x}{\| x\|} \]
is the desired measure.
(2) Since we have a measure $\mu$ supported by $B_{\Phi}$ with a standard
covariance matrix it is easy to see that
${\rm dim} ({\rm span}\{ B_{\Phi}\})\p =\p n$.
Consequently, the formal identity $\iota_{\Phi}$ is correctly defined and we
have $\| \iota^{-1}_{\Phi} \| \kl 1$ as well as $\pz \iota_{\Phi}^* ) \kl
\Phi (\iota_{\Phi}) \p = \p \sqrt{n}$ (see Lemma \ref{phi_p2}).
(3) For the first inequality we can use well-known volume estimates for the
corresponding $\ell$-norm, namely

\[  \frac{\ell (\iota_{\Phi})}{\sqrt{n}}
    \pl =\pl \kla \int_{S_{n-1}} \| x\|_{E_{\Phi}}^2 d\sigma_n (x)\mer^{1/2}
    \gl \kla \int_{S_{n-1}} \| x\|_{E_{\Phi}}^{-n} d\sigma_n (x) \mer^{-1/n}
    \pl = \pl  \kla \frac{{\rm vol}(\bzn)}{{\rm vol}(B_{\Phi})} \mer^{1/n}
    \pl ,\]

where $\sigma_n$ is the normalized Haar measure on the sphere $S_{n-1}$.
The second inequality of (3) follows from the results of Carl and Pajor,
see\cite{CP}(Corollary 1.4(b)), since

\[ \kla \frac{{\rm vol}(\bzn )}{{\rm vol}(B_{\Phi })} \mer^{1/n} \kl 2 \p
   e_n(\iota_{\Phi}) \pl,\]

where $e_n(\iota_{\Phi})$ denotes the $n$-th dyadic entropy number of
$\iota_{\Phi}$. \hfill $\Box$ \hz

We deduce the lower estimates for the $\Phi$-type and -cotype constants
which emphasize the special role of the spaces $\len$ and $\lin$ in
Theorem \ref{phi-type}.\hz

\begin{prop} \label{pos_results} Let $\fis\subset L_2^N$ be an orthonormal
system. Then, for some constant $\theta >0$ not depending on $\Phi$, $n$, and
$N$, the following holds true.
\begin{enumerate}
\item [(1)] There exists an $n$-dimensional subspace
      $E\subseteq \ell_{\infty}^N$ with
      \[ \fco E \pl \ge \pl \theta \pl \sqrt{\frac{\log (n+1)}{\log (N+1)}}
	 \sqrt{\frac{n}{\log (n+1)}} \pl.\]
\item [(2)] In the real situation there exists an $n$-dimensional quotient
      $\ell^N_1 /L$ with
      \[ \fty {\ell_1^N /L} \pl \ge \pl \theta \pl
	 \sqrt{\frac{1}{\log \kla \frac{N}{n} +1 \mer }} \pl \sqrt{n}\pl.\]

\end{enumerate}
\end{prop} \hz

$Proof:$ (1) We define the operator $R\in \com \ell_{\infty}^N)$ with
$\| R\| \le 1$ and $\sqrt{n} \kl \Phi(R)$ as in the proof of
Proposition \ref{smallcotype}. Then we consider the space $E\p:=\p R(\lzn)$
of dimension at most n. From $\ell(R) \kll c_0 \sqrt{\log (N+1)} \p \| R\|$
and the above inequality it follows that
\[ \sqrt{ \frac{n}{\log (n+1)} } \kl c_0\pl \fco E \pl
   \sqrt{ \frac{\log (N+1)}{\log(n+1)} } \pl.\]
(2) Since $L_2 (\Om ,p) = L_2^N$ we get
$B_{\Phi} = \overline{\rm absconv} \{ x_1 ,...,x_N\}$
for some $x_1,...,x_N \in S_{n-1}$ and $E_{\Phi} = \ell_1^N /L$.
Corollary 2.4(i) in \cite{CP} implies

\[ \kla \frac{{\rm vol}(B_{\Phi})}{{\rm vol}(\bzn )} \mer^{1/n} \le c\pl
   \sqrt{\frac{\log (\frac{N}{n} +1)}{n}} \]

for some absolute constant $c>0$. Using Lemma \ref{phi_to_mu} (3) we deduce

\for \sqrt{ \frac{n}{\ln \kla \frac{N}{n}+1 \mer }}
     \kl c \pl \frac{\ell(\iota_{\Phi})}{\sqrt{n}}
     \kl c \pl \fty {\efi} \pl \frac{ \Phi (\iota_{\Phi}) }
     {\sqrt{n}} \p = \p c \pl \fty {\efi}\pl.\\[-1.5cm] \mel \hfill $\Box$
\hz\hz

\begin{rem} Note that assertion (1) can formulated as follows in
the cotype situation:
As long as $N$ does not grow faster than a polynomial in $n$ there is an $n$
dimensional subspace of $\ell_{\infty}^N$ with "worst possible" cotype -$\Phi$
constant, although this constant can be quite small on $\lin$, see
Proposition \ref{smallcotype}. On the other hand assertion (2) in the type
situation
means: As long as $N \sim n$ there is an $n$ dimensional quotient of
$\ell_1^N$
with "worst possible" type - $\Phi$ constant, although this constant can be
bounded on $\ell_1$, see Proposition \ref{smalltype}.
\end{rem} \hz

\begin{rem} \label{ell_pif}
\begin{enumerate}
\item For any orthonormal system $\fis$ and any $u\in\com X)$ we have
      $\ell (u) \le \pif (u)$. This is known and follows
      (for example) from Lemma \ref{phi_to_mu} (1). In fact, we have
      \begin{eqnarray*}
      \pif (u)
      & = & \sqrt{n} \sup_{w\in {\cal O}_n} \kla \int_{\bzn} \| uwx\|^2
	     d\mu (x) \mer^{\frac{1}{2}}
      \ge \sqrt{n} \kla \int_{w\in {\cal O}_n} \int_{\bzn} \| uwx\|^2
	     d\mu (x) d w \mer^{\frac{1}{2}} \\
      & = & \sqrt{n} \kla \int_{\bzn} \| ux\|^2 d\sigma_n (x)
	    \mer^{\frac{1}{2}} ,
      \end{eqnarray*}
      where $dw$ stands for the integration with respect to the Haar measure
      on the group ${\cal O}(n)$ of the $n$-dimensional orthogonal matrices.
\item Let $\fis\subset\ell_2^N$ be an (real) orthonormal system such that
      \[ \|\summ_1^n\vare_k x_k\|_2\pl\le\pl c_0\pl\|\summ_1^n\phi_k
	 x_k\|_2 \pl\]
      holds for all Banach spaces $X$ and all $x_1,...,x_n \in X$. Then
      $N\ge n\kla e^{\frac{n}{(cc_0)^2}} -1\mer$ for some numerical
      constant $c>0$. To see this we consider $E_{\Phi} = \ell_1^N/L$ as
      in the proof of Proposition \ref{pos_results} and observe
      \[ \frac{r(\iota_{\Phi})}{\sqrt{n}} \gl \frac{1}{c_1} \pl \kla
	 \frac{{\rm vol}(\bzn)}{{\rm vol}(B_{\Phi})} \mer^{1/n}
	 \gl \frac{1}{c_1 c_2} \sqrt{\frac{n}{\log (N/n +1)}} \]
      according to Lemma \ref{phi_to_mu} (3) and \cite{CP} (Corollary 2.4).
      Since $\Phi (\iota_{\Phi}) = \sqrt{n}$ we get
      \[ \sqrt{n} = \Phi (\iota_{\Phi}) \gl \frac{r(\iota_{\Phi})}{c_0}
	 \gl \frac{1}{c_o c_1 c_2} \p \sqrt{n} \sqrt{\frac{n}{\log (N/n +1)}}
	 \]
      such that $\frac{n}{(c_0 c_1 c_2)^2} \kl \log \kla \frac{N}{n} +1
      \mer$.
\end{enumerate}
\end{rem} \hz

Finally, we are in a position to complete Remark \ref{nontrivial_ann}.\hz

\begin{prop} \label{nontrivial_sys} There exists an orthonormal system
$\Phi = (\phi_k)_1^{\infty} \subseteq (e^{ikt})_{k\in\nz}
\subseteq L_2(\Pi )$ such that
\[ \sup_n \frac{c_2 (\lin )}{c_{\Phi}(\lin )} \p = \p
   \sup_n \frac{t_2 (\len )}{t_{\Phi}(\len )} \p = \p \infty \hspace{.7cm}
   \mbox{and} \hspace{.7cm}
   c_{\Phi} (\ell_{\infty}) = t_{\Phi}(\ell_{\infty}) = \infty \p .\]
\end{prop}\hz

$Proof:$ We use the system constructed in \cite{BOUR}(Theorem 2). For
$2<q<\infty$ and $k=1,2,...$ there are chosen subsets
$S_k \subseteq \{ n: 2^k \le n < 2^{k+1} \}$ satisfying
$| S_k | = [4^{k/q}]$ such that for $\Lambda = \cup_{k=1}^{\infty} S_k$ the
system $\Phi = (e^{ikt})_{k\in \Lambda}$ is a $K_q$ system. Using
Lemma \ref{Kaq} we obtain

\[ \fco \lin \kl c_q \frac{n^{\frac{1}{q}}}{\sqrt{\log (n+1)}} K_q (\Phi)
   \pla\mbox{and}\pla
   \sup_n \frac{c_2 (\lin )}{c_{\Phi}(\lin )} \p = \infty . \]

Applying $t_{\Phi} (\len ) \kl K \sqrt{\log (n+1)} \fco \lin$ one also gets
$\sup_n \frac{t_2 (\len )}{t_{\Phi}(\len )} \p = \p \infty$ .
Now let us consider the system $\Psi_k := (e^{ilt})_{l\in S_k}$. It follows
from the Marcinkiewicz-Zygmund-inequality and a shift argument that

\[ \kla \frac{1}{2^k} \summ_{l=1}^{2^k} \noo \summ_{s+2^k-1\in S_k}
   e^{2\pi i \frac{sl}{2^k}} x_s \rrm^2 \mer^{\frac{1}{2}} \kl c
   \kla \int_0^{2\pi} \noo \summ_{s\in S_k} e^{ist} x_s \rrm^2
   \frac{dt}{2\pi} \mer^{\frac{1}{2}} \]

for all $x_1,...,x_{|S_k|} \in X$ and all Banach spaces $X$. Setting
$\Psi_k^0 := \kla \kla e^{2\pi i \frac{sl}{2^k}} \mer_{l=1}^{2^k}
\mer_{s=1}^{|S_k|} \subset L_2^{2^k}$ we get

\[ \pi_{\Psi_k^0} (T) \kl c \pi_{\Psi_k} (T) \hspace{1cm} \mbox{for all}
   \hspace{1cm} T\in\comxy .\]

Applying Proposition \ref{p2_few} we continue to

\[ \frac{1}{12} \kla \frac{|S_k|}{2^k} \mer^{\frac{1}{2}} \pi_2^{|S_k|} (T)
   \kl \pi_{\Psi_k^0} (T) \kl c \p \pi_{\Psi_k} (T) \]

for all $T\in\comxy$. Setting $n=2^k$ this reads as

\[ \frac{1}{12 c} \kla \frac{[n^{2/q}]}{n}\mer^{\frac{1}{2}} \p
   \pi_2^{[n^{2/q}]} (T) \kl \pi_{\Psi_k} (T)\p . \]

If $T$ is the embedding of $\ell_2^{[n^{2/q}]}$ into
$\ell_{\infty}^{[n^{2/q}]}$ this gives

\[ \frac{1}{12 c} \left [ n^{2/q} \right ]^{\frac{1}{2} +\frac{1}{2}}
   \p n^{-\frac{1}{2}} \kl
   \pi_{\Psi_k} \kla \iota_{2,\infty}^{[n^{2/q}]} \mer . \]

On the other hand

\[ \ell \kla \iota_{2,\infty}^{[n^{2/q}]} \mer \sim \sqrt{\log [n^{2/q}]} \]

which yields for $\frac{2}{q} -\frac{1}{2} >0$ ($q<4$) the equality
$\sup_n \fco \lin \p  = \p \infty$ . Finally, we have for any subsystem
$(\psi_k )_1^n \subset \Phi$ the estimates ($(e_k)$ is the standard basis of
$\lin$)

\[ \noo \summ_1^n \psi_k e_k \rrm_{L_2(\lin )} \p = \p 1
   \hspace{.7cm} \mbox{and} \hspace{.7cm}
   \noo \summ_1^n g_k e_k \rrm_{L_2(\lin )} \p \sim \p \sqrt{\log (n+1)} \]

such that $\sup_n t_{\Phi} (\lin ) = \infty$. \hfill $\Box$



\end{document}